\documentclass[11pt]{article}
\usepackage{bbm}
\usepackage{amsmath,amsthm,amssymb}
\usepackage{graphicx}
\usepackage{tikz}
\usepackage[utf8]{inputenc}
\usepackage{pgfplots}

\usepackage{setspace}
\usepackage{pict2e,color}
\usepackage{multirow}
\usepackage{booktabs}
\usepackage{amsfonts}
\usepackage{verbatim} 
\usepackage{sectsty}
\usepackage{url}
\usepackage{subfigure}
\usepackage{empheq}
\usepackage[normalem]{ulem}
\usepackage[titletoc,title]{appendix}
\usepackage[flushleft]{threeparttable}
\usepackage{soul} 
\usepackage{natbib}
\usepackage{algorithm,algorithmic}
\usepackage{etoolbox}
\preto\subequations{\ifhmode\unskip\fi}
\usepackage{booktabs}
\usepackage{multirow}
\usepackage{natbib}
\theoremstyle{definition}
\newtheorem{definition}{Definition}
\newtheorem{theorem}{Theorem}

\newtheorem{proposition}{Proposition}

\newtheorem{assumption}{Assumption}

\title{\Large Multistage Distributionally Robust Mixed-Integer Programming with Decision-Dependent Moment-Based Ambiguity Sets}

\begin{document}
\allowdisplaybreaks
\author{Xian Yu\thanks{Department of Industrial and Operations Engineering, University of Michigan at Ann Arbor, USA. Email: {\tt yuxian@umich.edu};}
    ~~~Siqian Shen\thanks{Corresponding author; Department of Industrial and Operations Engineering, University of Michigan at Ann Arbor, USA. Email: {\tt siqian@umich.edu}.}}
\date{}

\maketitle

\begin{abstract}
 We study multistage distributionally robust mixed-integer programs under endogenous uncertainty, where the probability distribution of stage-wise uncertainty depends on the decisions made in previous stages. We first consider two ambiguity sets defined by decision-dependent bounds on the first and second moments of uncertain parameters and by mean and covariance matrix that exactly match decision-dependent empirical ones, respectively. For both sets, we show that the subproblem in each stage can be recast as a mixed-integer linear program (MILP). Moreover, we extend the {general moment-based} ambiguity set in \citep{delage2010distributionally} to the multistage decision-dependent setting, and derive mixed-integer semidefinite programming (MISDP) reformulations of {stage-wise} subproblems. We develop methods for attaining lower and upper bounds of the optimal objective value of the multistage MISDPs, and approximate them using a series of MILPs. We deploy the Stochastic Dual Dynamic integer Programming (SDDiP) method for solving {the problem under the three ambiguity sets} with risk-neutral or risk-averse objective functions, and conduct numerical studies on {multistage facility-location instances having diverse sizes under different parameter and uncertainty settings. Our results show that the SDDiP quickly finds optimal solutions for moderate-sized instances under the first two ambiguity sets, and also finds good approximate bounds for the multistage MISDPs derived under the third ambiguity set. We also demonstrate the efficacy of incorporating decision-dependent distributional ambiguity in multistage decision-making processes.} 
\end{abstract}
   


\section{Introduction}
\label{sec:intro}
Data uncertainty appears ubiquitously in decision-making processes in practice, where system design and operational decisions {are made sequentially and} dynamically over a finite time horizon, to be adaptive to varying parameters (e.g., random customer demand, stochastic travel time). When using stochastic programming approaches, the goal is to optimize a certain measure of a random outcome (e.g., the expected cost of service operations)  given a fully known distribution of uncertain parameter. We refer to, e.g., \citet{birge2011introduction, shapiro2009lectures}, for detailed discussions about applications, formulations, and solution algorithms used in two-stage and multistage stochastic programming. On the other hand, robust optimization \citep{ben2009robust, bertsimas2011theory} provides an alternative way to make conservative decisions, and assumes that values of uncertain parameter may vary in a given constrained set, called ``uncertainty set.'' The resultant model seeks a solution that is feasible for any realization in the uncertainty set and optimal for the worst-case objective function. 

Recently, an approach that bridges the gap between robust optimization and stochastic programming is proposed to handle decision-making problems with ambiguously known distributions of uncertain parameter, namely, the distributionally robust optimization (DRO) approach. In DRO, optimal solutions are sought for the worst-case probability distribution within a family of candidate distributions, called an ``ambiguity set.'' A seminal paper by \citet{delage2010distributionally} focused on ambiguity sets defined by mean and covariance matrix, where they proved that a distributionally robust convex program can be reformulated as a semidefinite program and solved in polynomial time for a wide range of objective functions. They also quantified the relationship between the amount of data and the choice of moment-based ambiguity set parameters for achieving certain levels of solution conservatism. Recent DRO literature demonstrates that the ways of constructing the ambiguity sets can base on (i) empirical moments and their nearby regions \citep[see, e.g.,][]{mehrotra2014cutting, wagner2008stochastic, zhang2018ambiguous, delage2010distributionally}, and (ii) statistical distances between a candidate
distribution and a reference distribution, such as norm-based distance \citep[see][]{jiang2018risk}), $\phi$-divergence \citep[see][]{jiang2016data}, and Wasserstein metric \citep[see, e.g.,][]{esfahani2018data, blanchet2019quantifying, gao2016distributionally}. In this paper, we focus on moment ambiguity sets and extend them to multistage decision-dependent uncertainty settings, which we elaborate later. 

\citet{bertsimas2018adaptive} studied adaptive DRO in a dynamic setting, where decisions are adapted to the uncertain outcomes through stages. They focused on a class of second-order conic representable ambiguity sets and transformed the adaptive DRO problem to a classical robust optimization problem following linear decision rules. \citet{goh2010distributionally} studied a linear optimization problem under uncertainty which has  expectation terms in the objective function and constraints. The authors developed a new nonanticipative decision rule, which was more flexible than the linear decision rule, to find DRO solutions.

In practice, system parameters and therefore their uncertain features could depend on decisions made previously. For example, customer demand in various types of service industries, especially new service or service launched in a new market, is random and hard to predict due to lack of prior data. Its probability distribution can be greatly dependent on locations of service centers or facilities.  For example, consider carsharing or bikesharing services offered in metropolitan areas. Normally, one would sign up as a member only if she can easily find available cars or bikes nearby her work/home locations \citep[see][]{kung2018approximation}.
This type of uncertainty is called endogenous uncertainty, which has been extensively studied in the literature of dynamic programming \citep[see, e.g.,][]{webster2012approximate}), stochastic programming \citep[see, e.g.,][]{goel2006class, jonsbraaten1998class, lee2012newsvendor}) and robust optimization \citep[see, e.g.,][]{poss2013robust, spacey2012robust, hu2019data, lappas2017use, lappas2018robust, nohadani2018optimization}). Among them, \citet{webster2012approximate} proposed an approximate dynamic programming approach to solve a multistage global climate policy problem under decision-dependent uncertainties. \citet{goel2006class} studied a class of stochastic programs with decision-dependent parameters and presented a hybrid mixed-integer disjunctive programming formulation for these programs. \citet{poss2013robust} investigated robust combinatorial optimization with variable budgeted uncertainty, where the uncertain parameters belong to the image of multifunctions of the problem variables. They proposed a mixed-integer linear program (MILP) to reformulate the problem. Furthermore, \citet{vayanos2011decision} considered the process of revealing uncertain information being affected by previously made decisions, and proposed decision rules for stochastic programs with decision-dependent information discovery processes. \citet{vayanos2020robust} extended their methods to a robust optimization setting and performed numerical studies on instances of the active preference elicitation problem, solved for designing city security and crime control policies.

We consider multistage mixed-integer DRO models under endogenous uncertainty, of which the ambiguity sets are moment based and depend on previous stages' decisions. {The following papers also incorporate decision-dependent uncertainty into DRO formulations, but do not consider multistage, dynamic, nested formulations as the ones we will introduce in Sections \ref{sec:dro} and \ref{sec:ambiguity_3}.} \citet{noyan2018distributionally} considered a DRO problem, where the ambiguity sets are balls centered at a decision-dependent probability distribution. The measure they used is based on a class of earth mover's distances, 
including both total variation distance and Wasserstein metrics. Their models are nonconvex nonlinear programs, which are computationally intractable, and the authors specified several problem settings under which it is possible to obtain tractable formulations. {They demonstrated the results by solving small instances of a distributionally robust job scheduling problem that only involves 1 machine, 2 jobs, and 2 scenarios in the finite support of uncertain job-processing time.}
{\citet{luo2020distributionally} studied two-stage DRO models with decision-dependent ambiguity sets constructed using bounds on moments, covariance matrix, Wasserstein metric, Phi-divergence and Kolmogorov–Smirnov test. For the finite support case, they provide a small numerical example of a newsvendor problem where both the decision variable and uncertainty are 1-dimensional.} Recently, \citet{basciftci2019distributionally} considered a two-stage distributionally robust facility location problem, where mean and variance of the demand depend on the first-stage facility-opening decisions. The authors derived an equivalent MILP based on special problem structures and developed valid inequalities to improve the solution time when testing larger-sized instances {(with up to 10 facility locations, 20 demand sites, and 100 possible realizations in the support of demand)}. 

Regarding algorithms for multistage stochastic programs, \citet{pereira1991multi} were the first to develop the Stochastic Dual Dynamic Programming (SDDP) algorithm for efficiently computing multistage stochastic linear programs based on scenario tree representation of the dynamically realized uncertainty. We also refer the interested readers to \citet{philpott2008convergence, girardeau2014convergence, guigues2016convergence} for studies on the convergence of the SDDP algorithm under different problem settings. Recently, \citet{philpott2018distributionally} studied a variant of SDDP with a distributionally robust objective, where the ambiguity set is a Euclidean neighborhood of the nominal probability distribution. The authors showed its almost-sure convergence under standard assumptions and applied it to New Zealand hydrothermal electricity system. Stochastic Dual Dynamic integer Programming (SDDiP), firstly proposed by \citet{zou2019stochastic}, is an extension of SDDP to handle the nonconvexity arising in multistage stochastic integer programs. The essential differences are the new reformulations of subproblems in each stage and a new class of cuts derived for handling the integer variables. 

In this paper, we deploy risk-neutral expectation and risk-averse coherent-risk measures to interpret the objective functions in multistage DRO models with {decision-dependent endogenous uncertain parameter}. We consider three types of moment-based ambiguity sets respectively involving: Type 1 decision-dependent bounds on moments {(extended from one case of ambiguity sets in \citep{luo2020distributionally} for two-stage decision-dependent DRO models)}; Type 2 the mean vector and covariance matrix exactly matching decision-dependent empirical ones 
{(extended from the ambiguity set proposed by \citet{wagner2008stochastic} for general DRO models)}; and Type 3 the mean vector of uncertain parameters lying in an ellipsoid centered at a decision-dependent estimate mean vector, and the centered second-moment matrix lying in a positive semidefinite (psd) cone {(extended from the general moment ambiguity set in \citep{delage2010distributionally})}. For Type 1 and Type 2 ambiguity sets, we reformulate the problem as multistage stochastic MILPs, and for Type 3, we reformulate it as a multistage stochastic mixed-integer semidefinite program (MISDP). We then apply variants of the SDDiP approach for solving these reformulations or deriving objective bounds.

The main contributions of the paper are threefold. First, to our best knowledge, this paper is the first that handles mixed-integer DRO models under endogenous uncertainty in a multistage setting and derives reformulations that can be solved by off-the-shelf solvers. 
Second, the reformulation for Type 3 ambiguity set is a multistage MISDP, which cannot be optimized directly by any state-of-the-art integer-programming solvers. We derive both lower- and upper-bounds via Lagrangian relaxation and inner approximation, respectively, and numerically show that these bounds can approximate the optimal objective of the multistage problem well by having 4\% optimality gap in most instances given demand with high variation. Third, we successfully implement the SDDiP algorithm for handling both risk-neutral and risk-averse models and numerically evaluate the efficacy of our reformulations and bounds via testing diverse-sized problems (in terms of number of decision variables, constraints, stages in SDDiP and the support size). 

The remaining of this paper is organized as follows. {In Section \ref{sec:n-dddr}, we set up the formulation of a risk-neutral multistage decision-dependent DRO model with mixed-integer variables in each stage, and describe our problem assumptions.} In Section \ref{sec:dro}, we develop {exact MILP} reformulations and SDDiP algorithms for the {multistage decision-dependent DRO} models under Type 1 and Type 2 ambiguity sets. {In Section \ref{sec:ambiguity_3}, we develop MISDP reformulations and bounds for approximating the optimal objective for Type 3 ambiguity set. In Section \ref{sec:numerical}, we consider multistage facility-location instances having location-dependent demand and a finite set of periods for locating facilities. We demonstrate the finite convergence of the SDDiP algorithm, and present numerical results for instances with different sizes and parameter settings.} In Section \ref{sec:con}, we conclude the paper and state future research directions. 

Furthermore, we present all reformulations for the continuous support case in Appendix \ref{appen:continuous}, analysis of the risk-averse models under the three ambiguity sets in Appendix \ref{sec:dro_risk}, and details of all proofs in Appendix \ref{appen:allproofs}.

Throughout the paper, we use the following notation: The bold symbol will be used to denote a vector/matrix; for $n\in \mathbb{Z}_{+}$, the set $\{1,\ldots,n\}$ is represented by $[n]$; the Frobenius inner product $\text{trace} (\boldsymbol{A}^{\mathsf T} \boldsymbol{B})$ is denoted by $\boldsymbol{A}\bullet \boldsymbol{B}$.

\section{Problem Formulation and Assumptions}
\label{sec:n-dddr}
In the main paper, we focus on risk-neutral multistage decision-dependent distributionally robust mixed-integer programming models. (Due to similar analysis and results, we describe reformulations for the risk-averse models having coherent-risk-based objectives in Appendix \ref{sec:dro_risk}.) 

Consider a generic formulation of a multistage DRO problem with endogenous uncertainty and risk-neutral objectives as 
\begin{align}
\mbox{{\bf N-DDDR:}} & ~\nonumber \\
\min_{(\boldsymbol{x}_1,\boldsymbol{y}_1)\in X_1}&\Big\{g_1(\boldsymbol{x}_1,\boldsymbol{y}_1)+\max_{P_2\in \mathcal{P}_2(\boldsymbol{x}_1)}\mathbb{E}_{P_2}\Big[\min_{(\boldsymbol{x}_2,\boldsymbol{y}_2)\in X_2(\boldsymbol{x}_1,\boldsymbol{\xi}_2)}g_2(\boldsymbol{x}_2,\boldsymbol{y}_2)+\cdots\nonumber\\
&+\max_{P_{t}\in\mathcal{P}_{t}(\boldsymbol{x}_{t-1})}\mathbb{E}_{P_{t}}\Big[\min_{(\boldsymbol{x}_{t},\boldsymbol{y}_{t})\in X_{t}(\boldsymbol{x}_{t-1},\boldsymbol{\xi}_{t})}g_{t}(\boldsymbol{x}_{t},\boldsymbol{y}_{t})+\cdots\nonumber\\
&+\max_{P_T\in\mathcal{P}_T(\boldsymbol{x}_{T-1})}\mathbb{E}_{P_T}\Big[\min_{(\boldsymbol{x}_T,\boldsymbol{y}_T)\in X_T(\boldsymbol{x}_{T-1},\boldsymbol{\xi}_T)}g_T(\boldsymbol{x}_T,\boldsymbol{y}_T)\Big]\Big\}, \label{eq:1}
\end{align}
where $\boldsymbol{\xi}_t\in\Xi_t\subset\mathbb{R}^{J}$ is the random vector at stage $t$, for all $t=2,\ldots, T$. W.l.o.g., let $\Xi_1$ be a singleton, i.e., $\boldsymbol{\xi}_1$ is a deterministic vector. For $t > 1$, the probability of each uncertain parameter $\boldsymbol{\xi}_t$ is not known exactly, but lies
in an ambiguity set of probability distributions. Letting $\Xi=\Xi_1\times\Xi_2\times\cdots\times\Xi_T$, the evolution of $\boldsymbol{\xi}_t$ defines a probability space $(\Xi,\mathcal{F},P)$, and a filtration $\mathcal{F}_1\subset\mathcal{F}_2\subset\cdots\subset\mathcal{F}_T\subset\mathcal{F}$ such that each $\mathcal{F}_t$ corresponds to the information available up to (and including) the current stage $t$, with $\mathcal{F}_1=\{\emptyset,\Xi\}, \ \mathcal{F}_T=\mathcal{F}$.  We define binary state variable  $\boldsymbol{x}_t\in\{0,1\}^{I}$ to connect the consecutive two stages $t$ and $t+1$, and define integer/continuous stage variable $\boldsymbol{y}_t\in\mathbb{R}^{I\times J}$ which only appears at stage $t$. The feasible region for choosing decisions $(\boldsymbol{x}_t, \boldsymbol{y}_t)$ is $X_t(\boldsymbol{x}_{t-1},\boldsymbol{\xi}_t)\subset \{0,1\}^{I}\times\mathbb{R}^{I\times J}$, which depends on the values of decision $\boldsymbol{x}_{t-1}$ and random vector $\boldsymbol{\xi}_t$. Consider linear cost function  $g_t(\boldsymbol{x}_t,\boldsymbol{y}_t)$ and non-empty compact mixed-integer polyhedral feasible set $X_t(\boldsymbol{x}_{t-1},\boldsymbol{\xi}_t)$ for each $t\in[T]$. The ambiguity set at stage $t$ is denoted by $\mathcal{P}_t(\boldsymbol{x}_{t-1})$, which depends on the previous stage's decision variable $\boldsymbol{x}_{t-1}$, and $\mathcal{P}_t(\boldsymbol{x}_{t-1})\subset \mathcal{P}_t(\Xi_t,\mathcal{F}_t)$, denoting the set of probability distributions defined on $(\Xi_t,\mathcal{F}_t)$, for all $t=2,\ldots, T$.

The dynamic decision-making process is as follows:
\begin{align*}
&\underbrace{\text{decision}\ (\boldsymbol{x}_1, \boldsymbol{y}_1)}_{{\text{Stage $1$}}}\to \underbrace{{\text{worst-case}}\ (P_2) \to \text{observation}\ (\boldsymbol{\xi}_2)\to \text{decision}\ (\boldsymbol{x}_2, \boldsymbol{y}_2) }_{{\text{Stage $2$}}}\to
 \cdots\\ 
& \hspace{6ex}\to  \underbrace{{\text{worst-case}}\ (P_t) \to  \text{observation}\ (\boldsymbol{\xi}_{t})\to \text{decision}\ (\boldsymbol{x}_{t}, \boldsymbol{y}_t)}_{{\text{Stage $t$}}} \to \cdots \\
 &\hspace{18ex}\to  \underbrace{{\text{worst-case}}\ (P_T) \to  \text{observation}\ (\boldsymbol{\xi}_{T})\to \text{decision}\ (\boldsymbol{x}_{T}, \boldsymbol{y}_T)}_{{\text{Stage $T$}}}
\end{align*}
In the first stage, we make decisions $\boldsymbol{x}_1,\ \boldsymbol{y}_1$. The nature chooses the worst-case probability distribution $P_2\in\mathcal{P}_2(\boldsymbol{x}_1)$, under which the uncertain parameter $\boldsymbol{\xi}_2$ is observed and then make corresponding decisions $\boldsymbol{x}_2,\ \boldsymbol{y}_2$ in the second stage. This process continues until reaching stage $T$.

The Bellman equations for N-DDDR Model \eqref{eq:1} involve: 
\begin{equation*}
Q_1=\min_{(\boldsymbol{x}_1,\boldsymbol{y}_1)\in X_1} \left\{g_1(\boldsymbol{x}_1,\boldsymbol{y}_1)+\max_{P_2\in\mathcal{P}_2(\boldsymbol{x}_1)}\mathbb{E}_{P_2}[Q_2(\boldsymbol{x}_1,\boldsymbol{\xi}_2)]\right\},
\end{equation*}
\begin{equation}
Q_t(\boldsymbol{x}_{t-1},\boldsymbol{\xi}_t)=\min_{(\boldsymbol{x}_t,\boldsymbol{y}_t)\in X_t(\boldsymbol{x}_{t-1},\boldsymbol{\xi}_t)} \left\{g_t(\boldsymbol{x}_t,\boldsymbol{y}_t)+\max_{P_{t+1}\in\mathcal{P}_{t+1}(\boldsymbol{x}_t)}\mathbb{E}_{P_{t+1}}[Q_{t+1}(\boldsymbol{x}_t,\boldsymbol{\xi}_{t+1})]\right\},\label{eq:bell_neutral}
\end{equation}
for each $t=2,\ldots,T-1$, and
\begin{equation}
Q_T(\boldsymbol{x}_{T-1},\boldsymbol{\xi}_T)=\min_{(\boldsymbol{x}_T,\boldsymbol{y}_T)\in X_T(\boldsymbol{x}_{T-1},\boldsymbol{\xi}_T)} g_T(\boldsymbol{x}_T,\boldsymbol{y}_T).\nonumber
\end{equation}
Note that the Bellman equation in each stage $t\in[T-1]$ is a min-max problem. Therefore, our goal is to recast the inner maximization problem as a minimization problem and then reformulate the min-max model as a monolithic formulation.  Let $\hat{X}_t$ represent the feasible set $X_t$ projecting to the $\boldsymbol{x}_t$-space, i.e., $\boldsymbol{x}_t\in \hat{X}_t$ if and only if there exists $\boldsymbol{y}_t$ such that $(\boldsymbol{x}_t,\boldsymbol{y}_t)\in X_t$. We make the following assumptions in this paper. 

\begin{assumption}
The random vectors are stage-wise independent, i.e., $\boldsymbol{\xi}_t$ is stochastically independent of $\boldsymbol{\xi}_{[1,t-1]}=(\boldsymbol{\xi}_1,\ldots,\boldsymbol{\xi}_{t-1})^{\mathsf T}$, for all $t=2,\ldots,T$. 
\end{assumption}

\begin{assumption}
\label{assum:complete-recourse}
{The subproblem $Q_t(\boldsymbol{x}_{t-1},\boldsymbol{\xi}_t)$ in each stage $t$ is always feasible for any decision made in the constraint set $X_t$ and for every realization of the random vector $\boldsymbol{\xi}_t$ for all $t\in[T]$. That is, the problem has complete recourse.}
\end{assumption}

\begin{assumption}\label{assumption:finite}
For each $t=2,\ldots,T$, every probability distribution  $P_t\in\mathcal{P}_t(\boldsymbol{x}_{t-1})$ has a decision-independent support $\Xi_t:=\{\boldsymbol{\xi}_t^k\}_{k=1}^K$ with finite $K$ elements for all solution values $\boldsymbol{x}_{t-1}\in \hat{X}_{t-1}$. Each realization $\boldsymbol{\xi}^k_t$ is associated with a {decision-dependent} ambiguously known probability ${p_k(\boldsymbol{x}_{t-1})}$ satisfying $\sum_{k =1}^K {p_k(\boldsymbol{x}_{t-1})} = 1$. 
\end{assumption}

Assumption \ref{assum:complete-recourse} is for notation simplicity of the derivation and analysis of the SDDiP algorithm for solving reformulations of Model~\eqref{eq:1}. It is made w.l.o.g. as we can always penalize the violation of a certain constraint in the objective function by adding an additional penalty-related variable to the constraint. 

Assumption \ref{assumption:finite} is needed for deriving {efficient, finitely convergent} algorithms for multistage models. If we relax the assumption and allow continuous supports $\Xi_t, \ \forall t =2,\ldots,T$, the reformulations of N-DDDR under three ambiguity sets become semi-infinite programs with an infinite number of constraints and cannot be numerically tested. (We will present the corresponding reformulations in Theorems \ref{theorem:1-C}, \ref{theorem:2-C} and \ref{theorem:3-C} in Appendix \ref{appen:continuous}.) Therefore, we keep Assumption 3 in the main paper to derive reformulations of N-DDDR, and numerically evaluate their performance in Section \ref{sec:numerical}.

For notation simplicity, every discrete support $\Xi_t$ is assumed to have the same number of elements $K$ for $t=2,\ldots,T$. However, our model and solution approaches can be easily extended to settings with time-varying $K$. {Moreover, our setting can also accommodate the case of decision-dependent support with $\Xi_t(\boldsymbol{x}_{t-1})= \{\boldsymbol{\xi}_t^k\}_{k \in [K]\setminus \Omega}$, by letting a subset $\Omega \subset [K]$ of realizations to have zero probabilities, i.e., $p_k(\boldsymbol{x}_{t-1})=0,\ \forall k\in \Omega$, if our decision $\boldsymbol{x}_{t-1}$ will not lead to any of those specific realizations $\boldsymbol{\xi}^k_t, \forall k \in \Omega$ in stage $t$ for all $t=2,\ldots,T$.}

\section{Solving N-DDDR under Type 1 and Type 2 Ambiguity Sets}
\label{sec:dro}
We consider Types 1 and 2 ambiguity sets mentioned in Section \ref{sec:intro} for characterizing ambiguity sets $\mathcal P_2(\boldsymbol{x}_1),$ $\ldots,$ $\mathcal P_T(\boldsymbol{x}_{T-1})$, and will derive MILP reformulations and algorithms for {exactly optimizing} N-DDDR under these two ambiguity sets. 

\subsection{Reformulation under Type 1 Ambiguity Set}
\label{sec:ambiguity_1}
Following the settings of one ambiguity set studied by \citet{luo2020distributionally}, we bound all the moments by certain decision-dependent functions. In stage $t+1$, the random vector is $\boldsymbol{\xi}_{t+1}=(\xi_{t+1,1},\ldots,\xi_{t+1,J})^{\mathsf T}\in \mathbb{R}^{J}$ where $\xi_{t+1,j}$ represents the $j$-th uncertain parameter.  We consider $m$ different moment functions 
$\boldsymbol{f}:=(f_1(\boldsymbol{\xi}_{t+1}),\ldots,f_m(\boldsymbol{\xi}_{t+1}))^{\mathsf T}$. 
Then for each $s=1,\ldots, m$, 
$$f_s(\boldsymbol{\xi}_{t+1})=(\xi_{t+1,1})^{k_{s1}}(\xi_{t+1,2})^{k_{s2}}\cdots(\xi_{t+1,J})^{k_{sJ}},$$ 
where $k_{sj}$ is a non-negative integer indicating the power of $\xi_{t+1,j}$ for the $s$-th moment function. 
The lower and upper bounds are defined by $\boldsymbol{l}(\boldsymbol{x}_t) := (l_1(\boldsymbol{x}_t),\ldots, l_m(\boldsymbol{x}_t))^{\mathsf T}$ and $ \boldsymbol{u}(\boldsymbol{x}_t) := (u_1(\boldsymbol{x}_t),\ldots, u_m(\boldsymbol{x}_t))^{\mathsf T}$, respectively. For each $t\in[T-1]$, a discrete Type 1 ambiguity set with Assumption \ref{assumption:finite} is:  
\begin{align}\label{eq:ambi}
\mathcal{P}_{t+1}^{D_1}(\boldsymbol{x}_t):=\left\{\boldsymbol{p}\in\mathbb{R}^K\ |\ \underline{\boldsymbol{p}}(\boldsymbol{x}_t)\le \boldsymbol{p}\le \bar{\boldsymbol{p}}(\boldsymbol{x}_t),\ \boldsymbol{l}(\boldsymbol{x}_t)\le \sum_{k=1}^Kp_k\boldsymbol{f}(\boldsymbol{\xi}_{t+1})\le\boldsymbol{u}(\boldsymbol{x}_t)\right\},
\end{align}
where $\underline{\boldsymbol{p}}(\boldsymbol{x}_t)$ and $\bar{\boldsymbol{p}}(\boldsymbol{x}_t)$ are the given lower and upper bounds of the candidate true probability $\boldsymbol{p}$, which are decision-dependent. Following the derivations in \citep{luo2020distributionally} for reformulating a two-stage decision-dependent DRO model, we generalize their results for the multistage setting and reformulate Bellman equation \eqref{eq:bell_neutral} below in Theorem \ref{theorem:1}. 
Note that $\boldsymbol{p}$ can be ensured as a probability distribution by setting one of the moment functions $\boldsymbol{f}$, lower and upper bounds $\boldsymbol{l}$ and $\boldsymbol{u}$ to be 1 (which then enforces $\sum_{k=1}^K p_k = 1$). The details are given in equations \eqref{eq:normal} and \eqref{eq:ambi_ineq:a} later and without loss of generality, we do not include \eqref{eq:ambi_ineq:a} specifically in \eqref{eq:ambi}.
We also describe a continuous version of $ \mathcal{P}_{t+1}^{D_1}(\boldsymbol{x}_t)$ and the resulting reformulation in Appendix \ref{appen:continuous}.

\begin{theorem} 
\label{theorem:1}
If for any feasible $\boldsymbol{x}_t\in \hat{X}_t$, the ambiguity set defined in \eqref{eq:ambi} is non-empty, then the Bellman equation \eqref{eq:bell_neutral} can be reformulated as: 
\begin{subequations}\label{Bellman}
\begin{align}
Q_t(\boldsymbol{x}_{t-1},\boldsymbol{\xi}_t)=\min_{\boldsymbol{\alpha},\boldsymbol{\beta},{\boldsymbol{\underline{\gamma}},\boldsymbol{\bar{\gamma}},}\boldsymbol{x}_t,\boldsymbol{y}_t}\quad &g_t(\boldsymbol{x}_t,\boldsymbol{y}_t)-{\boldsymbol{\alpha}}^{\mathsf T} \boldsymbol{l}(\boldsymbol{x}_t)+{\boldsymbol{\beta}}^{\mathsf T} \boldsymbol{u}(\boldsymbol{x}_t) {-\underline{\boldsymbol{\gamma}}^{\mathsf T}\underline{\boldsymbol{p}}(\boldsymbol{x}_t) + \bar{\boldsymbol{\gamma}}^{\mathsf T}\bar{\boldsymbol{p}}(\boldsymbol{x}_t)}\label{eq:obj_bellman}\\
\text{s.t.}\quad&(-\boldsymbol{\alpha}+\boldsymbol{\beta})^{\mathsf T} \boldsymbol{f}(\boldsymbol{\xi}_{t+1}^k) {-\underline{\gamma}_k+\bar{\gamma}_k}\ge Q_{t+1}(\boldsymbol{x}_t,\boldsymbol{\xi}_{t+1}^k),\ \forall k\in[K], \label{eq:constr_bellman}\\
&(\boldsymbol{x}_t,\boldsymbol{y}_t)\in X_t(\boldsymbol{x}_{t-1},\boldsymbol{\xi}_t),\\
& \boldsymbol{\alpha},\ \boldsymbol{\beta},\  {\boldsymbol{\underline{\gamma}},\ \boldsymbol{\bar{\gamma}}}\ge 0.
\end{align}
\end{subequations}
\end{theorem}
 
The proof of Theorem \ref{theorem:1} is presented in Appendix \ref{appen:allproofs}.
Note that there exist nonlinear terms in both objective function \eqref{eq:obj_bellman} and constraints \eqref{eq:constr_bellman} (e.g., ${\boldsymbol{\alpha}}^{\mathsf T} \boldsymbol{l}(\boldsymbol{x}_t),\ {\boldsymbol{\beta}}^{\mathsf T} \boldsymbol{u}(\boldsymbol{x}_t)$) and we explore special structures of $\mathcal{P}^{D_1}_{t+1} (\boldsymbol{x}_{t})$ to speed up the computation. 
For the first and second moments of each parameter, we consider their lower and upper bounds as follows:
\small
\begin{subequations}
\begin{align}
&f_1(\boldsymbol{\xi}_{t+1})=1,\ l_1(\boldsymbol{x}_t)=u_1(\boldsymbol{x}_t)=1, \label{eq:normal} \\ 
&f_{1+j}(\boldsymbol{\xi}_{t+1})=\xi_{t+1,j},\ l_{1+j}(\boldsymbol{x}_t)=\mu_j(\boldsymbol{x}_t)-\epsilon_j^{\mu}, \
u_{1+j}(\boldsymbol{x}_t)=\mu_j(\boldsymbol{x}_t)+\epsilon_j^{\mu},\ \forall j\in[J] ,\label{eq:first_moment}\\
&f_{1+J+j}(\boldsymbol{\xi}_{t+1})=(\xi_{t+1,j})^2,\
l_{1+J+j}(\boldsymbol{x}_t)=S_j(\boldsymbol{x}_t)\underline{\epsilon}_j^S,\
u_{1+J+j}(\boldsymbol{x}_t)=S_j(\boldsymbol{x}_t)\bar{\epsilon}_j^S,\ \forall j\in[J]. \label{eq:second_moment}
\end{align}
\end{subequations}
\normalsize
Here, \eqref{eq:normal} is a normalization constraint to ensure that $P$ is a probability distribution. Equations \eqref{eq:first_moment} and \eqref{eq:second_moment} demonstrate the first and second moment functions for each parameter, respectively. When the first moment function is used, $l_{1+j}(\boldsymbol{x}_t)$ and $u_{1+j}(\boldsymbol{x}_t)$ bound the mean of parameter $\xi_{t+1,j}$ in an $\epsilon_j^{\mu}$-interval of the empirical mean function $\mu_j(\boldsymbol{x}_t)$ for all $j\in[J]$. Similarly, $l_{1+J+j}(\boldsymbol{x}_t)$ and $u_{1+J+j}(\boldsymbol{x}_t)$ bound the second moment of parameter $\xi_{t+1,j}$ via scaling the empirical second moment function $S_j(\boldsymbol{x}_t)$ for all $j\in[J]$. { In the rest of our analysis, we set $\underline{\boldsymbol{p}}(\boldsymbol{x}_t) = \boldsymbol{0},\ \bar{\boldsymbol{p}}(\boldsymbol{x}_t)=\boldsymbol{1}$ for any feasible $\boldsymbol{x}_t$, and focus on specially designed forms of $\mu_j(\boldsymbol{x}_t)$ and $S_j(\boldsymbol{x}_t)$ to derive a computable reformulation of Model \eqref{Bellman}. 
}
We first specify $2J+1$ constraints in the ambiguity set \eqref{eq:ambi} as:
\small
\begin{subequations}
\label{eq:ambi_ineq}
\begin{align}
\mathcal{P}^{D_1}_{t+1} (\boldsymbol{x}_{t})=\Biggl \{\boldsymbol{p}\in\mathbb{R}_{+}^K\ |&\sum_{k=1}^Kp_{k}=1, \label{eq:ambi_ineq:a}\\
& \mu_j(\boldsymbol{x}_t)-\epsilon_j^{\mu}\le\sum_{k=1}^K p_{k}\xi_{t+1,j}^k\le \mu_j(\boldsymbol{x}_t)+\epsilon_j^{\mu},\ \forall j\in[J],\\
& S_j(\boldsymbol{x}_t)\underline{\epsilon}_j^S\le\sum_{k=1}^K p_{k}(\xi_{t+1,j}^k)^2\le S_j(\boldsymbol{x}_t)\bar{\epsilon}_j^S,\ \forall j\in[J]\Biggr \}, 
\end{align}
\end{subequations}
\normalsize
where for each $j\in[J]$, the empirical first moment $\mu_j(\boldsymbol{x}_t)$ and second moment $S_j(\boldsymbol{x}_t)$ affinely depend on decisions $\boldsymbol{x}_t$, such that
\begin{align*}
&\mu_j(\boldsymbol{x}_t)=\bar\mu_j\left(1+\sum_{i=1}^I\lambda_{ji}^{\mu}x_{ti}\right),\\
&S_j(\boldsymbol{x}_t)=(\bar\mu_j^2+\bar\sigma_j^2)\left(1+\sum_{i=1}^I\lambda_{ji}^Sx_{ti}\right),
\end{align*}
where the empirical mean and standard deviation of the $j$-th uncertain parameter are denoted by $\bar\mu_j,\ \bar\sigma_j$, respectively. Here by assumption, the first and second moments will increase when any of the state variable $x_{ti}$ changes from 0 to 1. Parameters $\lambda_{ji}^{\mu},\ \lambda_{ji}^S\in\mathbb{R}_{+}$ respectively represent the degree about how $x_{ti}=1$ may affect the values of the first and second moments of $\xi_{t+1,j}$ for each $j\in[J]$. Following this assumption, the mean and variance of customer demand may increase if there are more facilities open nearby, and the respective increasing rates are measured by $\lambda^{\mu}_{ji}$ and $\lambda_{ji}^S$. Depending on specific applications and problem contexts, the values of $\lambda^\mu$'s and $\lambda^S$'s can be set differently. Also note that for notation simplicity, $\lambda^u_{ji}$ and $\lambda^S_{ji}$ are the same for all stages $t\in[T]$. Our models and approaches can also accommodate time-varying $\lambda^\mu$- or $\lambda^S$-values. 

We further rewrite the recursive function $Q_{t+1}(\boldsymbol{x}_t,\boldsymbol{\xi}_{t+1}^k)$ as $Q_{t+1}^k$ for notation simplicity. Using the ambiguity set defined in \eqref{eq:ambi_ineq}, the Bellman equation \eqref{Bellman} becomes 
\small
\begin{subequations}\label{eq:bell_ineq}
\begin{align}
Q_t(\boldsymbol{x}_{t-1},\boldsymbol{\xi}_t)=\min_{\boldsymbol{\alpha},\boldsymbol{\beta},\boldsymbol{x}_t,\boldsymbol{y}_t}\quad &g_t(\boldsymbol{x}_t,\boldsymbol{y}_t)-\alpha_{1}-\sum_{j=1}^J\alpha_{2j}(\bar\mu_j-\epsilon_j^{\mu})-\sum_{j=1}^J\sum_{i=1}^I\lambda_{ji}^{\mu}\bar\mu_j\alpha_{2j}x_{ti}\nonumber\\
\quad& -\sum_{j=1}^J\alpha_{3j}(\bar\mu_j^2+\bar\sigma_j^2)\underline{\epsilon}^S_j-\sum_{j=1}^J\sum_{i=1}^I\lambda^S_{ji}\underline{\epsilon}^S_j(\bar\mu_j^2+\bar\sigma_j^2)\alpha_{3j}x_{ti}\nonumber\\
\quad&+\beta_{1}+\sum_{j=1}^J\beta_{2j}(\bar\mu_j+\epsilon_j^{\mu})+\sum_{j=1}^J\sum_{i=1}^I\lambda_{ji}^{\mu}\bar\mu_j\beta_{2j}x_{ti}\nonumber\\
\quad& +\sum_{j=1}^J\beta_{3j}(\bar\mu_j^2+\bar\sigma_j^2)\bar{\epsilon}^S_j+\sum_{j=1}^J\sum_{i=1}^I\lambda^S_{ji}\bar{\epsilon}^S_j(\bar\mu_j^2+\bar\sigma_j^2)\beta_{3j}x_{ti}\label{eq:Bellman2}\\
\text{s.t.}\quad&-\alpha_{1}+\beta_{1}+\sum_{j\in [J]}\xi_{t+1,j}^k(-\alpha_{2j}+\beta_{2j})+\sum_{j\in[J]}(\xi_{t+1,j}^k)^2(-\alpha_{3j}+\beta_{3j})\ge Q_{t+1}^k,\nonumber\\
& \hspace{12ex}\forall k\in[K],\label{eq:value_1}\\
&(\boldsymbol{x}_t,\boldsymbol{y}_t)\in X_t(\boldsymbol{x}_{t-1},\boldsymbol{\xi}_t),\label{eq:xy}\\
& \boldsymbol{\alpha},\ \boldsymbol{\beta}\ge 0.\nonumber
\end{align}
\end{subequations}
\normalsize
Given binary valued $x_{ti}$,  we provide exact reformulations of the bilinear terms $z^{\alpha_2}_{tji}=\alpha_{2j}x_{ti},\ z^{\alpha_3}_{tji}=\alpha_{3j}x_{ti}, \ z^{\beta_2}_{tji}=\beta_{2j}x_{ti}, \ z^{\beta_3}_{tji}=\beta_{3j}x_{ti}$ in objective \eqref{eq:Bellman2} using McCormick envelopes $M^{\alpha_2}_{tji},\ M^{\alpha_3}_{tji}, \ M^{\beta_2}_{tji},\ M^{\beta_3}_{tji}$ for all $i\in[I],\ j\in[J]$. {(We omit constraint details of all the McCormick envelopes here and also in the remaining reformulations as they follow standard procedures, which can be found in, e.g., \citet{mccormick1976computability}.)} 

Then, following the multi-cut version of SDDiP algorithm \citep{zou2019stochastic}, at iteration $\ell$, we replace the value function $Q_{t+1}^k$ by under-approximation cuts:
\begin{equation}
\theta_{t}^k\ge v_{t+1}^{lk}+(\boldsymbol{\pi}_{t+1}^{lk})^{\mathsf T}\boldsymbol{x}_t,\ \forall k\in[K],\ l\in[\ell],\label{eq:cuts}
\end{equation}
where cut coefficients $\{(v_{t+1}^{lk}, \boldsymbol{\pi}_{t+1}^{lk})\}_{k=1}^K$ are evaluated at stage $t+1$ in the backward step at each iteration $l \in [\ell]$ with $\boldsymbol{\pi}_{t+1}^{lk}$ being the optimal solution to a Lagrangian dual problem of model \eqref{eq:bell_ineq} and $v_{t+1}^{lk}=\mathcal{L}_{t+1}^k(\boldsymbol{\pi}_{t+1}^{lk})$ being the value of the Lagrangian dual function. Then we obtain an under-approximation of the Bellman equation \eqref{eq:bell_ineq} as
\small
\begin{align}
\underline{Q}_t(\boldsymbol{x}_{t-1},\boldsymbol{\xi}_t)=\min_{\substack{\boldsymbol{\alpha},\boldsymbol{\beta},\boldsymbol{x}_t,\boldsymbol{y}_t,\boldsymbol{\theta}_t\\\boldsymbol{z}^{\alpha_2}, \boldsymbol{z}^{\alpha_3}, \boldsymbol{z}^{\beta_2},\boldsymbol{z}^{\beta_3}}}\quad &g_t(\boldsymbol{x}_t,\boldsymbol{y}_t)-\alpha_{1}-\sum_{j=1}^J\alpha_{2j}(\bar\mu_j-\epsilon_j^{\mu})-\sum_{j=1}^J\sum_{i=1}^I\lambda_{ji}^{\mu}\bar\mu_jz^{\alpha_2}_{tji}\nonumber\\
\quad& -\sum_{j=1}^J\alpha_{3j}(\bar\mu_j^2+\bar\sigma_j^2)\underline{\epsilon}^S_j-\sum_{j=1}^J\sum_{i=1}^I\lambda_{ji}\underline{\epsilon}^S_j(\bar\mu_j^2+\bar\sigma_j^2)z^{\alpha_3}_{tji}\nonumber\\
\quad&+\beta_{1}+\sum_{j=1}^J\beta_{2j}(\bar\mu_j+\epsilon_j^{\mu})+\sum_{j=1}^J\sum_{i=1}^I\lambda_{ji}^{\mu}\bar\mu_jz^{\beta_2}_{tji}\nonumber\\
\quad& +\sum_{j=1}^J\beta_{3j}(\bar\mu_j^2+\bar\sigma_j^2)\bar{\epsilon}^S_j+\sum_{j=1}^J\sum_{i=1}^I\lambda_{ji}\bar{\epsilon}^S_j(\bar\mu_j^2+\bar\sigma_j^2)z^{\beta_3}_{tji}\label{eq:milp}\\
\text{s.t.}\quad&\mbox{\eqref{eq:xy},\ \eqref{eq:cuts}}, \nonumber\\
&-\alpha_{1}+\beta_{1}+\sum_{j\in[J]}\xi_{t+1,j}^k(-\alpha_{2j}+\beta_{2j})+\sum_{j\in[J]}(\xi_{t+1,j}^k)^2(-\alpha_{3j}+\beta_{3j})\ge \theta_{t}^k,\nonumber\\
& \hspace{12ex}\forall k\in[K],\nonumber\\
&(z^{\alpha_2}_{tji},\alpha_{2j},x_{ti})\in M^{\alpha_2}_{tji},\ \forall i\in[I],\ j\in[J],\nonumber\\
&(z^{\alpha_3}_{tji},\alpha_{3j},x_{ti})\in M^{\alpha_3}_{tji},\ \forall i\in[I],\ j\in[J],\nonumber\\
&(z^{\beta_2}_{tji},\beta_{2j},x_{ti})\in M^{\beta_2}_{tji},\ \forall i\in[I],\ j\in[J],\nonumber\\
&(z^{\beta_3}_{tji},\beta_{3j},x_{ti})\in M^{\beta_3}_{tji},\ \forall i\in[I],\ j\in[J],\nonumber\\
& \boldsymbol{\alpha},\ \boldsymbol{\beta}\ge 0. \nonumber
\end{align}
\normalsize
The above under-approximation \eqref{eq:milp} is an MILP. Therefore, we can apply SDDiP using Lagrangian cuts to optimize the original N-DDDR model \eqref{eq:1} with its stage-wise subproblem reformulations \eqref{eq:milp}, given Type 1 ambiguity set.

\subsection{Reformulation under Type 2 Ambiguity Set}
\label{sec:dro_match}
In the previous section, we consider Type 1 ambiguity set defined by decision-dependent bounds on each moment separately, whereas in reality, there may be correlations between different parameters. In this case, we rely on estimates of the true mean and covariance matrix and consider ambiguity sets defined by matching empirical mean $\boldsymbol{\mu}(\boldsymbol{x}_t)\in\mathbb{R}^{J}$ and covariance matrix $\boldsymbol{\Sigma}(\boldsymbol{x}_t)\in\mathbb{R}^{J\times J}$ exactly. For each $t\in[T-1]$, we consider Type 2 ambiguity set having a discrete support of uncertain parameter, given by
\small
\begin{subequations}
\label{eq:ambi_match}
  \begin{align}
 \mathcal{P}^{D_2}_{t+1}(\boldsymbol{x}_t):=\Biggl \{\boldsymbol{p}\in \mathbb{R}^k\ |\  & \sum_{k=1}^K p_k =1, \\
&\sum_{k=1}^Kp_k\boldsymbol{\xi}_{t+1}^k=\boldsymbol{\mu}(\boldsymbol{x}_t),\\
& \sum_{k=1}^Kp_k(\boldsymbol{\xi}_{t+1}^k-\boldsymbol{\mu}(\boldsymbol{x}_t))(\boldsymbol{\xi}_{t+1}^k-\boldsymbol{\mu}(\boldsymbol{x}_t))^{\mathsf T}=\boldsymbol{\Sigma}(\boldsymbol{x}_t)\Biggr \}. 
\end{align}
\end{subequations}
\normalsize
Theorem \ref{theorem:2} demonstrates a reformulation of Bellman equation \eqref{eq:bell_neutral} given Type 2 ambiguity set \eqref{eq:ambi_match}. 
\begin{theorem}
\label{theorem:2}
If for any feasible $\boldsymbol{x}_t\in \hat{X}_t$, the ambiguity set defined in \eqref{eq:ambi_match} is non-empty, then the Bellman equation \eqref{eq:bell_neutral} can be reformulated as
\small
\begin{subequations}\label{Bellman_2}
 \begin{align}
Q_t(\boldsymbol{x}_{t-1},\boldsymbol{\xi}_t)=\min_{\boldsymbol{x}_t,\boldsymbol{y}_t,s,\boldsymbol{u},\boldsymbol{Y}}\quad &g_t(\boldsymbol{x}_t,\boldsymbol{y}_t)+s+\boldsymbol{u}^{\mathsf T}\boldsymbol{\mu}(\boldsymbol{x}_t)+\boldsymbol{\Sigma}(\boldsymbol{x}_t)\bullet \boldsymbol{Y} \label{eq:4}\\
\text{s.t.}\quad&s+\boldsymbol{u}^{\mathsf T}\boldsymbol{\xi}_{t+1}^k +(\boldsymbol{\xi}_{t+1}^k-\boldsymbol{\mu}(\boldsymbol{x}_t))(\boldsymbol{\xi}_{t+1}^k-\boldsymbol{\mu}(\boldsymbol{x}_t))^{\mathsf T}\bullet \boldsymbol{Y}\ge Q_{t+1}^k,\nonumber\\
&\hspace{12ex}\forall k\in[K],\\
&(\boldsymbol{x}_t,\boldsymbol{y}_t)\in X_t(\boldsymbol{x}_{t-1},\boldsymbol{\xi}_t).\nonumber
\end{align}
\end{subequations}
\normalsize
\end{theorem}

A detailed proof of Theorem \ref{theorem:2} is presented in Appendix \ref{appen:allproofs}, in which we apply strong duality to recast the inner maximization problem in \eqref{eq:bell_neutral} as a minimization problem and combine it with the outer minimization problem. Furthermore, assume that the elements in $\boldsymbol{\mu}(\boldsymbol{x}_t),\ \boldsymbol{\Sigma}(\boldsymbol{x}_t)$ are affine in $\boldsymbol{x}_t$, i.e.,
\begin{subequations}
\begin{align}
&\mu_j(\boldsymbol{x}_t)=\bar\mu_j(1+\sum_{i=1}^I\lambda_{ji}^{\mu}x_{ti}),\ \forall j\in[J]\label{eq:linear1}\\
&\boldsymbol{\Sigma}(\boldsymbol{x}_t)=\bar{\boldsymbol{\Sigma}}(1+\sum_{i=1}^I{\lambda}_{i}^{cov}x_{ti}),\label{eq:linear2}
\end{align}
\end{subequations}
where $\bar{\boldsymbol{\mu}}$ is the nominal mean vector and $\bar{\boldsymbol{\Sigma}}$ is a psd matrix representing the nominal covariance matrix. Then 
\begin{subequations}
\begin{align}
& \boldsymbol{u}^{\mathsf T}\boldsymbol{\mu}(\boldsymbol{x}_t) = \sum_{j=1}^J \bar{\mu}_ju_j(1+\sum_{i=1}^I \lambda_{ji}^{\mu}x_{ti}),\label{eq:bilinear1}\\
&\boldsymbol{\Sigma}(\boldsymbol{x}_t)\bullet Y=\sum_{j=1}^J\sum_{j^{\prime}=1}^J \bar{\Sigma}_{jj^{\prime}}(1+\sum_{i=1}^I{\lambda}_{i}^{cov}x_{ti})Y_{jj^{\prime}},\ \label{eq:bilinear}\\
&\boldsymbol{\mu}(\boldsymbol{x}_t)\boldsymbol{\mu}(\boldsymbol{x}_t)^{\mathsf T}\bullet Y=\sum_{j=1}^J\sum_{j^{\prime}=1}^J\bar{\mu}_j\bar\mu_{j^{\prime}}(1+\sum_{i=1}^I \lambda_{ji}^{\mu} x_{ti})(1+\sum_{i^{\prime}=1}^I \lambda_{j^{\prime}i^{\prime}}^{\mu} x_{ti^{\prime}})Y_{jj^{\prime}}. \ \label{eq:trilinear}
\end{align}
\end{subequations}
Both \eqref{eq:bilinear1} and \eqref{eq:bilinear} contain bilinear terms and \eqref{eq:trilinear} contains trilinear terms. Since $x_{ti}, \ \forall i\in[I]$ are binary variables, we can provide exact reformulations of the bilinear terms $w_{tij}=x_{ti}u_j,\ z_{tijj^{\prime}}=x_{ti}Y_{jj^{\prime}}$, and trilinear terms $v_{tii^{\prime}jj^{\prime}}=x_{ti}x_{ti^{\prime}}Y_{jj^{\prime}}$ for all $t\in[T],\ i,i^{\prime}\in[I],\ j,j^{\prime}\in[J]$ using McCormick envelopes $M^w_{tij}, \ M^z_{tijj^{\prime}}, \  M^v_{tii^{\prime}jj^{\prime}}$. 
Applying the same cutting planes in \eqref{eq:cuts}, we obtain an under-approximation $\underline{Q}_t(\boldsymbol{x}_{t-1},\boldsymbol{\xi}_t)$ of the Bellman equation \eqref{Bellman_2} as: 
\small \begin{align*}
\min_{\substack{\boldsymbol{x}_t,\boldsymbol{y}_t,s,\boldsymbol{u},\boldsymbol{Y}\\\boldsymbol{w},\boldsymbol{z},\boldsymbol{v}}}\quad &g_t(\boldsymbol{x}_t,\boldsymbol{y}_t)+s+\sum_{j=1}^J \bar{\mu}_ju_j+\sum_{i=1}^I\sum_{j=1}^J \bar{\mu}_j\lambda_{ji}^{\mu}w_{tij}+\sum_{j=1}^J\sum_{j^{\prime}=1}^J \bar{\Sigma}_{jj^{\prime}}Y_{jj^{\prime}}+\sum_{i=1}^I\sum_{j=1}^J\sum_{j^{\prime}=1}^J \bar{\Sigma}_{jj^{\prime}}{\lambda}_{i}^{cov}z_{tijj^{\prime}} \\
\text{s.t.}\quad& \mbox{\eqref{eq:xy},\ \eqref{eq:cuts}}\nonumber\\
&s+\boldsymbol{u}^{\mathsf T}\boldsymbol{\xi}_{t+1}^k + \boldsymbol{\xi}_{t+1}^k (\boldsymbol{\xi}_{t+1}^k)^{\mathsf T}\bullet \boldsymbol{Y} - \sum_{j=1}^J\sum_{j^{\prime}=1}^J \xi_{t+1,j}^k \bar{\mu}_{j^{\prime}}\Big(Y_{jj^{\prime}}+Y_{j^{\prime}j}+\sum_{i=1}^I \lambda_{j^{\prime}i}^{\mu}(z_{tij^{\prime}j}+z_{tijj^{\prime}})\Big)\\
&+\sum_{j=1}^J\sum_{j^{\prime}=1}^J\bar{\mu}_j\bar{\mu}_{j^{\prime}}\Big(Y_{jj^{\prime}}+\sum_{i=1}^I(\lambda_{j^{\prime}i}^{\mu}+\lambda_{ji}^{\mu})z_{tijj^{\prime}}+\sum_{i=1}^I\sum_{i^{\prime}=1}^I\lambda_{ji}^{\mu}\lambda_{j^{\prime}i^{\prime}}^{\mu}v_{tii^{\prime}jj^{\prime}}\Big)\ge \theta_{t}^k,\ \forall k\in[K],\\
& (w_{tij},x_{ti},u_j)\in M^w_{tij},\ \forall i\in[I],\ j\in[J],\\
& (z_{tijj^{\prime}},x_{ti},Y_{ij})\in M^z_{tijj^{\prime}},\ \forall i\in[I],\ j,j^{\prime}\in[J],\\
& (v_{tii^{\prime}jj^{\prime}},x_{ti^{\prime}},z_{tijj^{\prime}})\in M^v_{tii^{\prime}jj^{\prime}},\ \forall i,i^{\prime}\in[I],\ j,j^{\prime}\in[J], 
\end{align*}\normalsize
which is an MILP and we can again deploy the SDDiP approach for {optimally solving the N-DDDR model \eqref{eq:1}}.

\section{Solving N-DDDR under Type 3 Ambiguity Set}
\label{sec:ambiguity_3}

Now we focus on the general moment-based ambiguity set for decision-dependent DRO models, and derive reformulations and algorithms for N-DDDR under Type 3 ambiguity set, where the mean vector of uncertain parameters lies in an ellipsoid centered at an affinely decision-dependent estimate mean vector, and the second-moment matrix lies in a psd cone defined by an affinely decision-dependent matrix. Specifically, for all $t\in [T-1]$, letting $\gamma, \ \eta$ be coefficients controlling the size of the ambiguity set, we have
\small
\begin{subequations}
\label{eq:ambi_bound}
 \begin{align}
\mathcal{P}_{t+1}^{D_3}(\boldsymbol{x}_t):=\Biggl \{\boldsymbol{p}\in\mathbb{R}^K\ | &\sum_{k=1}^K p_k=1,\label{eq:ambi_bound1}\\ 
&\left(\sum_{k=1}^K p_k\boldsymbol{\xi}_{t+1}^k-\boldsymbol{\mu}(\boldsymbol{x}_t)\right)^{\mathsf T}\boldsymbol{\Sigma}(\boldsymbol{x}_t)^{-1}\left(\sum_{k=1}^K p_k\boldsymbol{\xi}_{t+1}^k-\boldsymbol{\mu}(\boldsymbol{x}_t)\right) \le \gamma,\\
& \sum_{k=1}^Kp_k(\boldsymbol{\xi}_{t+1}^k-\boldsymbol{\mu}(\boldsymbol{x}_t))(\boldsymbol{\xi}_{t+1}^k-\boldsymbol{\mu}(\boldsymbol{x}_t))^{\mathsf T}\preceq \eta \boldsymbol{\Sigma}(\boldsymbol{x}_t)\Biggr \}.\label{eq:ambi_bound4}
\end{align}
\end{subequations}
\normalsize

\subsection{Mixed-integer Semidefinite Programming Reformulation} 
Theorem \ref{theorem:3} describes a reformulation of Bellman equation \eqref{eq:bell_neutral} under Type 3 ambiguity set \eqref{eq:ambi_bound}.

\begin{theorem}
\label{theorem:3}
Suppose that the Slater's constraint qualification conditions
are satisfied, i.e., for any feasible $\boldsymbol{x}_t\in \hat{X}_t$, there exists a vector $\boldsymbol{p} = (p_1,p_2,\ldots,p_K)^{\mathsf T}$ such that $\sum_{k=1}^Kp_{k}=1$, $\left(\sum_{k=1}^K p_k\boldsymbol{\xi}_{t+1}^k-\boldsymbol{\mu}(\boldsymbol{x}_t)\right)^{\mathsf T}\boldsymbol{\Sigma}(\boldsymbol{x}_t)^{-1}\left(\sum_{k=1}^K p_k\boldsymbol{\xi}_{t+1}^k-\boldsymbol{\mu}(\boldsymbol{x}_t)\right)< \gamma$, and $\sum_{k=1}^Kp_k(\boldsymbol{\xi}_{t+1}^k-\boldsymbol{\mu}(\boldsymbol{x}_t))(\boldsymbol{\xi}_{t+1}^k-\boldsymbol{\mu}(\boldsymbol{x}_t))^{\mathsf T}\prec \eta \boldsymbol{\Sigma}(\boldsymbol{x}_t)$.
Using the ambiguity set defined in \eqref{eq:ambi_bound}, the Bellman equation \eqref{eq:bell_neutral} can be recast as
\small
\begin{subequations}\label{eq:bell_bound1}
 \begin{align}
Q_t(\boldsymbol{x}_{t-1},\boldsymbol{\xi}_t)=\min_{\boldsymbol{x}_t,\boldsymbol{y}_t,s,\boldsymbol{Z},\boldsymbol{Y}}\quad &g_t(\boldsymbol{x}_t,\boldsymbol{y}_t)+s+\boldsymbol{\Sigma}(\boldsymbol{x}_t)\bullet \boldsymbol{z}_1 -2\boldsymbol{\mu}(\boldsymbol{x}_t)^{\mathsf T}\boldsymbol{z}_2+\gamma z_3 +\eta\boldsymbol{\Sigma}(\boldsymbol{x}_t)\bullet \boldsymbol{Y}\\
\text{s.t.}\quad&s-2\boldsymbol{z}_2^{\mathsf T}\boldsymbol{\xi}_{t+1}^k+(\boldsymbol{\xi}_{t+1}^k-\boldsymbol{\mu}(\boldsymbol{x}_t))(\boldsymbol{\xi}_{t+1}^k-\boldsymbol{\mu}(\boldsymbol{x}_t))^{\mathsf T}\bullet Y\ge Q_{t+1}^k,\ \forall k\in[K],\\
& \boldsymbol{Z}=\begin{pmatrix} \boldsymbol{z}_1&\boldsymbol{z}_2\\\boldsymbol{z}_2^{\mathsf T}&z_3\end{pmatrix}\succeq 0,\ \boldsymbol{Y}\succeq 0,\\
&(\boldsymbol{x}_t,\boldsymbol{y}_t)\in X_t(\boldsymbol{x}_{t-1},\boldsymbol{\xi}_t).
\end{align}
\end{subequations}
\normalsize
\end{theorem}

A detailed proof of Theorem \ref{theorem:3} is given in Appendix \ref{appen:allproofs}. The key idea is to use the Lagrangian function and apply strong duality to recast the inner maximization problem in \eqref{eq:bell_neutral} as a minimization problem. We still assume the linear dependence of $\boldsymbol{\mu}(\boldsymbol{x}_t), \boldsymbol{\Sigma}(\boldsymbol{x}_t)$ on $\boldsymbol{x}_t$, as shown in \eqref{eq:linear1} and \eqref{eq:linear2}. Because $x_{ti}, \ i\in[I]$ are binary variables, we can provide exact reformulations of the bilinear terms $w_{tijj^{\prime}}=x_{ti}z_{1,jj^{\prime}},\ u_{tij}=x_{ti}z_{2j},\ R_{tijj^{\prime}}=x_{ti}Y_{jj^{\prime}}$, and trilinear terms $v_{tii^{\prime}jj^{\prime}}=x_{ti}x_{ti^{\prime}}Y_{jj^{\prime}}$  using McCormick envelopes $M^w_{tijj^{\prime}},\ M^u_{tij},\ M^R_{tijj^{\prime}},\ M^v_{tii^{\prime}jj^{\prime}}$ for all $t\in[T],\ i,i^{\prime}\in[I],\ j,j^{\prime}\in[J]$. 
Overall, the Bellman equation \eqref{eq:bell_bound1} can be recast as $Q_t(\boldsymbol{x}_{t-1},\boldsymbol{\xi}_t)=$
\small
\begin{subequations}\label{eq:bell_sdp}
\begin{align}
\min_{\substack{\boldsymbol{x}_t,\boldsymbol{y}_t,s,\boldsymbol{Z},\boldsymbol{Y}\\\boldsymbol{w},\boldsymbol{u},\boldsymbol{R},\boldsymbol{v}}} \quad&g_t(\boldsymbol{x}_t,\boldsymbol{y}_t)+s+\sum_{j=1}^J\sum_{j^{\prime}=1}^J\bar{\Sigma}_{jj^{\prime}}z_{1,jj^{\prime}}+\sum_{j=1}^J\sum_{j^{\prime}=1}^J \sum_{i=1}^I \bar{\Sigma}_{jj^{\prime}}{\lambda}_{i}^{cov}w_{tijj^{\prime}}\nonumber\\
\quad&-2\Big(\sum_{j=1}^J \bar{\mu}_jz_{2j}+\sum_{i=1}^I\sum_{j=1}^J\bar{\mu}_j\lambda_{ji}^{\mu}u_{tij}\Big)  +\gamma z_3+\eta\big(\sum_{j=1}^J\sum_{j^{\prime}=1}^J \bar{\Sigma}_{jj^{\prime}}Y_{jj^{\prime}}+\sum_{j=1}^J\sum_{j^{\prime}=1}^J\sum_{i=1}^I \bar{\Sigma}_{jj^{\prime}}{\lambda}_{i}^{cov}R_{tijj^{\prime}}\Big)\nonumber\\
\text{s.t.}\quad&s-2\boldsymbol{z}_2^{\mathsf T}\boldsymbol{\xi}_{t+1}^k + \boldsymbol{\xi}_{t+1}^k (\boldsymbol{\xi}_{t+1}^k)^{\mathsf T}\bullet \boldsymbol{Y} - \sum_{j=1}^J\sum_{j^{\prime}=1}^J \xi_{t+1,j}^k \bar{\mu}_{j^{\prime}}\Big(Y_{jj^{\prime}}+Y_{j^{\prime}j}+\sum_{i=1}^I \lambda_{j^{\prime}i}^{\mu}(R_{tij^{\prime}j}+R_{tijj^{\prime}})\Big)\nonumber\\
&+\sum_{j=1}^J\sum_{j^{\prime}=1}^J\bar{\mu}_j\bar{\mu}_{j^{\prime}}\Big(Y_{jj^{\prime}}+\sum_{i=1}^I(\lambda_{ji}^{\mu}+\lambda_{j^{\prime}i}^{\mu})R_{tijj^{\prime}}+\sum_{i=1}^I\sum_{i^{\prime}=1}^I\lambda_{ji}^{\mu}\lambda_{j^{\prime}i^{\prime}}^{\mu}v_{tii^{\prime}jj^{\prime}}\Big)\ge Q_{t+1}^k,\ \forall k\in[K],\label{eq:value}\\
&(\boldsymbol{x}_t,\boldsymbol{y}_t)\in X_t(\boldsymbol{x}_{t-1},\boldsymbol{\xi}_t),\label{eq:x}\\
& (w_{tijj^{\prime}},x_{ti},z_{1,jj^{\prime}})\in M^w_{tijj^{\prime}},\ \forall i\in[I],\ j,j^{\prime}\in[J],\label{eq:mc1}\\
& (u_{tij},x_{ti},z_{2j})\in M^u_{tij},\ \forall i\in[I],\ j\in[J],\\
& (R_{tijj^{\prime}},x_{ti},Y_{jj^{\prime}})\in M^R_{tijj^{\prime}},\ \forall i\in[I],\ j,j^{\prime}\in[J],\\
& (v_{tii^{\prime}jj^{\prime}},x_{ti^{\prime}},R_{tijj^{\prime}})\in M^v_{tii^{\prime}jj^{\prime}},\ \forall i,i^{\prime}\in[I],\ j,j^{\prime}\in[J],\label{eq:mc4}\\
& \boldsymbol{Z}=\begin{pmatrix} \boldsymbol{z}_1&\boldsymbol{z}_2\\\boldsymbol{z}_2^{\mathsf T}&z_3\end{pmatrix}\succeq 0,\ \boldsymbol{Y}\succeq 0.
\end{align}
\end{subequations}
\normalsize
For notation simplicity, we rewrite the linear objective function as $\tilde{g}_t(\boldsymbol{x}_t,\boldsymbol{y}_t, s,\boldsymbol{Z},\boldsymbol{Y},\boldsymbol{w},\boldsymbol{u},\boldsymbol{R})$ and the linear function on the left-hand side of Constraint \eqref{eq:value} as $f_t( s,\boldsymbol{Z},\boldsymbol{Y},\boldsymbol{R},\boldsymbol{v},\boldsymbol{\xi}_{t+1}^k)$. We fold all linear constraints \eqref{eq:x}--\eqref{eq:mc4} into set $\tilde{X}_t$. Then model \eqref{eq:bell_sdp} becomes:
\begin{subequations}
\label{eq:sdp}
\begin{align}
Q_t(\boldsymbol{x}_{t-1},\boldsymbol{\xi}_t)=\min_{\substack{\boldsymbol{x}_t,\boldsymbol{y}_t,s,\boldsymbol{Z},\boldsymbol{Y}\\\boldsymbol{w},\boldsymbol{u},\boldsymbol{R},\boldsymbol{v}}}\quad &\tilde{g}_t(\boldsymbol{x}_t,\boldsymbol{y}_t, s,\boldsymbol{Z},\boldsymbol{Y},\boldsymbol{w},\boldsymbol{u},\boldsymbol{R})\\
\text{s.t.}\quad&f_t( s,\boldsymbol{Z},\boldsymbol{Y},\boldsymbol{R},\boldsymbol{v},\boldsymbol{\xi}_{t+1}^k)\ge Q_{t+1}^k,\ \forall k\in[K],\label{eq:value2}\\
&(\boldsymbol{x}_t,\boldsymbol{y}_t,\boldsymbol{Z},\boldsymbol{Y},\boldsymbol{w},\boldsymbol{u},\boldsymbol{R}, \boldsymbol{v})\in \tilde{X}_t(\boldsymbol{x}_{t-1},\boldsymbol{\xi}_t),\\
& \boldsymbol{Z}\succeq 0 ,\ \boldsymbol{Y}\succeq 0.\label{eq:z}
\end{align}
\end{subequations}

\subsection{Derivation and Computation of Bounds for Multistage MISDPs} 
To solve \eqref{eq:sdp}, we aim to replace the value function $Q_{t+1}(\boldsymbol{x}_t, \boldsymbol{\xi}_{t+1}^k)$ in \eqref{eq:value2} by some under-approximation linear cuts, which will result in a multistage stochastic MISDP. The MISDP itself is difficult to solve directly due to the nature of semidefinite programs with integer variables. To our best knowledge, no solvers can directly optimize MISDP. For example, BNB and CUTSDP are two internal mixed-integer conic programming solvers in YALMIP \citep{lofberg2004yalmip}, which rely on relaxing integrality/semidefinite cones during iterative processes but not solve them exactly. If we want to leverage SDDiP with Lagrangian cuts, an MILP is needed in each stage. In the next two subsections, two methods are proposed to tackle this issue. In Section \ref{sec:relax}, we solve a Lagrangian relaxation, which provides valid cuts and the procedures will produce a lower bound on the optimal objective value of the original multistage problem. In Section \ref{sec:upper}, we approach the problem by inner approximating MISDPs via MILPs so that we can apply SDDiP with Lagrangian cuts directly on the resultant multistage MILP. The gaps of these two approaches are demonstrated numerically in Section \ref{sec:numerical}, to show the efficacy of the bounds.

\subsubsection{Lower bounding via Relaxed Lagrangian Cuts}
\label{sec:relax}
In the forward step, we solve the MISDPs \eqref{eq:sdp} for all stages $t\in[T-1]$ with current approximations of the value functions. Then in the backward step, at iteration $\ell$ of stage $t$, our goal is to find under-approximation linear cuts with coefficients $\{(v_t^{\ell k},\pi_t^{\ell k})\}_{k=1}^K$ for value function $Q_{t}(\boldsymbol{x}_{t-1}, \boldsymbol{\xi}_{t}^k)$ such that $Q_{t}(\boldsymbol{x}_{t-1},\boldsymbol{\xi}_{t}^k)\ge v_{t}^{\ell k} + (\boldsymbol{\pi}_{t}^{\ell k})^{\mathsf T}\boldsymbol{x}_{t-1}$ for all $\boldsymbol{x}_{t-1}\in\{0,1\}^I$. 
Following \citep{zou2019stochastic}, we make a copy of the state variable $\boldsymbol{z}_t=\boldsymbol{x}_{t-1}$ and then relax it to get a Lagrangian function. Specifically, at iteration $\ell$, for each realization $\boldsymbol{\xi}_t^{k}$, we solve the following relaxation problem in the backward step:
\begin{align*}
\mathcal{L}_t^{k}(\boldsymbol{\pi}_t)=\min_{\substack{\boldsymbol{x}_t,\boldsymbol{y}_t,\boldsymbol{z}_t,s,\boldsymbol{Z},\boldsymbol{Y}\\\boldsymbol{\theta}_{t},\boldsymbol{w},\boldsymbol{u},\boldsymbol{R},\boldsymbol{v}}}\quad &\tilde{g}_t(\boldsymbol{x}_t,\boldsymbol{y}_t, s,\boldsymbol{Z},\boldsymbol{Y},\boldsymbol{w},\boldsymbol{u},\boldsymbol{R}) - \boldsymbol{\pi}_t^{\mathsf T}\boldsymbol{z}_t\\
\text{s.t.}\quad&f_t( s,\boldsymbol{Z},\boldsymbol{Y},\boldsymbol{R},\boldsymbol{v},\boldsymbol{\xi}_{t+1}^{k^{\prime}})\ge \theta_t^{k^{\prime}},\ \forall k^{\prime}\in[K],\\
& \mbox{\eqref{eq:cuts}}\nonumber\\
&(\boldsymbol{x}_t,\boldsymbol{y}_t,\boldsymbol{Z},\boldsymbol{Y},\boldsymbol{w},\boldsymbol{u},\boldsymbol{R}, \boldsymbol{v})\in \tilde{X}_t(\boldsymbol{z}_{t},\boldsymbol{\xi}_t^{k}),\\
& \boldsymbol{Z}\succeq 0,\ \boldsymbol{Y}\succeq 0.
\end{align*}
A collection of cuts given by the coefficients $\{(v_t^{\ell k},\pi_t^{\ell k})\}_{k=1}^K$ is generated, where $\pi_t^{\ell k}\in\mathbb{R}^I$ is any real vector and $v_t^{\ell k}=\mathcal{L}_t^{k}(\pi_t^{\ell k})$. We name this collection of cuts the \textit{Relaxed Lagrangian Cuts} because it does not require the coefficient $\boldsymbol{\pi}$ to be the optimal solution to the Lagrangian dual problem.

\begin{proposition}
\label{prop:1}
The collection of Relaxed Lagrangian Cuts $\{(v_t^{lk},\pi_t^{lk})\}_{k=1}^K$ is valid because the true value function is bounded from below by these cuts for all $\boldsymbol{x}_{t-1}$, i.e., $Q_t(\boldsymbol{x}_{t-1},\boldsymbol{\xi}_t^k)\ge v_t^{lk} + (\boldsymbol{\pi}_t^{lk})^{\mathsf T}\boldsymbol{x}_{t-1}$ for all $\boldsymbol{x}_{t-1}\in\{0,1\}^I$.
\end{proposition}
The proof is similar to the one of Theorem 3 in \citep{zou2019stochastic} and it is omitted here. 

As a result, SDDiP algorithm with Relaxed Lagrangian Cuts provides a lower bound on the original multistage stochastic MISDP. However, because the Relaxed Lagrangian Cuts are not necessarily tight, our algorithm is not guaranteed to converge to an optimal solution. {In Section \ref{sec:numerical}, the tightness of the bounds is verified numerically based on diverse instances with different problem sizes and parameter settings.} 

\subsubsection{Upper bounding via inner approximating MISDP by MILPs}
\label{sec:upper}

We also propose to inner approximate psd cones by polyhedrons to obtain valid upper bounds for the MISDPs \eqref{eq:sdp}.
\begin{definition}
A symmetric matrix A is diagonally dominant (dd) if $a_{ii}\ge \sum_{j\not =i}|a_{ij}|$ for all $i$.
\end{definition}
We can further define a set of cones parameterized by a matrix $\boldsymbol{U}\in\mathbb{R}^{n\times n}$:
\begin{align*}
DD(\boldsymbol{U}):=\{\boldsymbol{M}\in S_n\ | \ \boldsymbol{M}=\boldsymbol{U}^{\mathsf T}\boldsymbol{Q}\boldsymbol{U} \ \text{for some dd matrix} \ \boldsymbol{Q}\},
\end{align*}
where $S_n$ represents the set of real symmetric $n\times n$ matrices. Optimizing over $DD(\boldsymbol{U})$ is a linear program since $\boldsymbol{U}$ is fixed and the associated constraints are linear in $\boldsymbol{M}$ and $\boldsymbol{Q}$. Moreover, the matrices in $DD(\boldsymbol{U})$ are all psd, i.e., $\forall \boldsymbol{U},\ DD(\boldsymbol{U})\subset P_n$, where $P_n$ represents the set of $n\times n$ psd matrices.

Then, following similar ideas in \citep{ahmadi2017sum}, one natural way is to replace the conditions $\boldsymbol{Z}\succeq 0,\ \boldsymbol{Y} \succeq 0$ by $\boldsymbol{Z}\in DD(\boldsymbol{U}),\ \boldsymbol{Y}\in DD(\boldsymbol{V})$ for some fixed matrices $\boldsymbol{U},\boldsymbol{V}$ in the forward step. This will provide us an upper bound on the value function $Q_t(\boldsymbol{x}_{t-1},\boldsymbol{\xi}_t)$, given by
\begin{align*}
\overline{Q}_t(\boldsymbol{x}_{t-1},\boldsymbol{\xi}_t)=\min_{\substack{\boldsymbol{x}_t,\boldsymbol{y}_t,s,\boldsymbol{Z},\boldsymbol{Y}\\\boldsymbol{w},\boldsymbol{u},\boldsymbol{R},\boldsymbol{v}}}\quad &\tilde{g}_t(\boldsymbol{x}_t,\boldsymbol{y}_t, s,\boldsymbol{Z},\boldsymbol{Y},\boldsymbol{w},\boldsymbol{u},\boldsymbol{R})\\
\text{s.t.}\quad&f_t( s,\boldsymbol{Z},\boldsymbol{Y},\boldsymbol{R},\boldsymbol{v},\boldsymbol{\xi}_{t+1}^k)\ge \overline{Q}_{t+1}(\boldsymbol{x}_t, \boldsymbol{\xi}_{t+1}^k),\ \forall k\in[K],\\
&(\boldsymbol{x}_t,\boldsymbol{y}_t,\boldsymbol{Z},\boldsymbol{Y},\boldsymbol{w},\boldsymbol{u},\boldsymbol{R}, \boldsymbol{v})\in \tilde{X}_t(\boldsymbol{x}_{t-1},\boldsymbol{\xi}_t),\\
& \boldsymbol{Z}\in DD(\boldsymbol{U}),\ \boldsymbol{Y}\in DD(\boldsymbol{V}).
\end{align*}
Then in the backward step, we can construct the Lagrangian cuts on the stage-wise MILPs. As a result, the optimal objective value of the resultant multistage MILP will serve as an upper bound of the original multistage MISDP.

In Appendix \ref{sec:dro_risk}, we generalize the risk-neutral objective functions in the N-DDDR model \eqref{eq:1} to risk-averse ones based on coherent risk measures. We present reformulations of the risk-averse multistage decision-dependent DRO problems under Types 1, 2, 3 ambiguity sets and derive SDDiP algorithms or bounds, similar to the results in Sections \ref{sec:dro} and \ref{sec:ambiguity_3}.

\section{Numerical Studies}
\label{sec:numerical}
We use instances of a multistage facility-location problem \citep[see, e.g.,][]{xian2019value} for validating our reformulations and algorithms. In these instances, consider $1,\ldots, I$ potential facilities and $1,\ldots, J$ customer sites. We define binary decision variable $x_{ti}$, $\forall t\in [T], \ i \in [I]$, such that $x_{ti}=1$ if a facility is open at location $i$ in stage $t$, and $x_{ti}=0$ otherwise. Decision variable $y_{tij}$ represents the flow of products from facility $i$ to customer site $j$ in stage $t$. The random vector at stage $t$ is $\boldsymbol{\xi}_t=(\xi_{t1},\ldots,\xi_{tJ})^{\mathsf T}$, representing the demand in each customer site at stage $t$. {Then, in model N-DDDR \eqref{eq:1}, the objective function at stage $t$ is defined as $g_t(\boldsymbol{x}_t, \boldsymbol{y}_t)= \sum_{i=1}^I\sum_{j=1}^J c_{ij}y_{tij}-\sum_{j=1}^J R_j\sum_{i=1}^I y_{tij}$, where it minimizes the total transportation cost minus the total revenue, and $c_{ij},\ R_j$ denote the unit transportation cost from facility $i$ to customer site $j$ and revenue for meeting one unit demand at customer site $j$, for all $i\in[I],\ j\in[J]$, respectively. The stage-wise feasibility set $X_{t}(\boldsymbol{x}_{t-1},\boldsymbol{\xi}_{t})$ for each $t\in[T]$ consists of the following constraints:
\small
\begin{subequations}
\label{model:deter}
\begin{align}
\quad& \sum_{i=1}^I y_{tij} \le \xi_{tj}, \quad\forall j\in[J], \label{eq:demand}\\
&	\sum_{j=1}^J y_{tij}\le h_{ti}\sum_{\tau=1}^t x_{\tau i},\quad\forall i\in[I], \label{eq:capacity}\\
&   \sum_{i=1}^I f_{ti}(x_{ti}-x_{t-1,i})\le N, \label{eq:budget}\\
&    x_{ti}\ge x_{t-1,i}, \quad\forall i\in[I], \label{eq:re}
\\
&	x_{ti}\in \{0,1\}, \quad\forall i\in[I],\\
&   y_{tij}\in \mathbb{Z}_{+}, \quad\forall i\in[I],\ j\in[J].
\end{align}
\end{subequations}
\normalsize
where \eqref{eq:demand} and \eqref{eq:capacity} require that the total shipment to a customer site/from a facility in each stage cannot exceed the demand/capacity of that customer site/facility, respectively. Constraints \eqref{eq:budget} imply that the building cost in each stage cannot exceed a given budget $N$, and according to \eqref{eq:re}, any open facilities cannot be removed.}

In all our tests, we randomly sample $I$ potential facilities and $J$ customer sites on a $100\times 100$ grid. The transportation costs between facilities and customer sites are calculated by their Manhattan distances divided by 4, i.e., $c_{ij}=\textrm{dist}(i,j)/4,\ \forall i\in[I],\ j\in[J]$. We set the building costs $f_{ti}=100$ for all $i\in[I]$ and $t\in[T]$. In each stage $t$, we set budget $N=100$, and all the facilities have the same capacity $h_{ti}=1000$ for all $t\in[T],\ i\in[I]$. The revenue for meeting one unit demand is set to $R_j=100$ for all $j \in[J]$.
  The empirical demand mean $\bar{\mu}_j$ is drawn uniformly between 20 and 40 for each $j\in[J]$, and the empirical standard deviation $\bar{\sigma}_j$ is set to $\bar{\mu}_j\times \bar{\rho}$, where we vary the coefficient $\bar{\rho}$ to represent different demand variations later. Then, for the uncertain demand $\xi_{tj}$, we sample $K$ data points following $\mathcal{N}(\bar{\mu}_j,\bar{\sigma}_j^2)$ for all $j\in[J]$, to construct the discrete support $\Xi_t$ in each stage $t\in[T]$.
  
In Section \ref{appen:exact}, we test small instances with $T=2$ stages, $I=3$ facilities and $J=1$ or $2$ customer site(s) for each of the three ambiguity sets. Specifically, we compare using the SDDiP algorithm for solving each reformulation of the N-DDDR model with an algorithm that enumerates all feasible first-stage solutions and optimizes the corresponding second-stage DRO models to seek optimal solutions. We show that both the optimal solutions and objective values of these two approaches are the same under the first two ambiguity sets, confirming the finite convergence of SDDiP algorithm to the true optimum. 

In Sections \ref{sec:compu_type1} and \ref{sec:compu_type3}, we test the SDDiP algorithm for solving reformulations given by Type 1 and Type 3 ambiguity sets, respectively,  on larger-sized instances by increasing values of $T$, $I$, $J$ and parameters used in SDDiP.

Our experiments utilize YALMIP toolbox in MATLAB \citep{lofberg2004yalmip} for modeling, where MOSEK is used to directly solve the stage-wise MILPs, and CUTSDP is used to solve MISDPs. All numerical experiments are conducted on a Windows 2012 Server with 128 GB RAM and an Intel 2.2 GHz processor.

\subsection{Results of Small Instances and Finite Convergence of SDDiP}
\label{appen:exact}

\subsubsection{Results of Type 1 ambiguity set on two-stage instances}
We first consider N-DDDR model with $T=2$ stages, $I=3$ facilities and $J=1$ customer site. For Type 1 ambiguity set \eqref{eq:ambi_ineq} in Section \ref{sec:ambiguity_1},  we set the empirical first and second moments as $\bar{\boldsymbol{\mu}}=10,\ \bar{\boldsymbol{\sigma}}=0.1$, and the bounding parameters as $\epsilon^{\mu}=5,\ \underline{\epsilon}^S=0.5,\ \bar{\epsilon}^S=1.5$. We evaluate four different patterns with fixed $\lambda^\mu$- and $\lambda^S$-values given in Table \ref{tab:pattern_1}.
\begin{table}[ht!]
    \centering
        \caption{Patterns with varying $\lambda^\mu$-/$\lambda^S$-values that interpret how decisions affect mean/variance.}
    \begin{tabular}{c|r|r}
    \hline
        Pattern & $(\lambda_{11}^{\boldsymbol{\mu}},\lambda_{12}^{\boldsymbol{\mu}},\lambda_{13}^{\boldsymbol{\mu}})^{\mathsf T}$ & $(\lambda_{11}^{\boldsymbol{\sigma}},\lambda_{12}^{\boldsymbol{\sigma}},\lambda_{13}^{\boldsymbol{\sigma}})^{\mathsf T}$ \\
        \hline
        1-1 & $(0.9,0.5,0.1)^{\mathsf T}$ & $(0.5,0.5,0.5)^{\mathsf T}$\\
        1-2 & $(0.5,0.5,0.5)^{\mathsf T}$ & $(0.9,0.5,0.1)^{\mathsf T}$\\
        1-3 & $(0.1,0.1,0.1)^{\mathsf T}$ & $(0.9,0.5,0.1)^{\mathsf T}$\\
        1-4 & $(0.5,0.9,0.1)^{\mathsf T}$ & $(0.9,0.5,0.1)^{\mathsf T}$\\
        \hline
    \end{tabular}
    \label{tab:pattern_1}
\end{table}
For each pattern, we first solve the two-stage min-max  formulation of N-DDDR by enumerating on all feasible first-stage solutions and each second-stage problem is directly optimized by MOSEK solver. We then apply SDDiP algorithm to solve both the N-DDDR and the decision-independent counterpart (N-DIDR) with all $\lambda^\mu$- and $\lambda^S$-values set to 0, where the algorithm iteratively builds cuts to approximate the first-stage value function.
Table \ref{tab:type1-2} demonstrates the performance of the above three models under different patterns. Each column under ``Two-stage enumeration'' displays the cost with the corresponding first-stage $\boldsymbol{x}$-solution ($\boldsymbol{x}$-sol.), where we mark the optimal solution in bold. The rest of the columns record the optimal objective values and optimal solutions of the N-DDDR and N-DIDR models, respectively. 

\begin{table}[ht!]
  \centering
  \caption{Results of different models using Type 1 ambiguity set}
  \resizebox{\textwidth}{!}{
    \begin{tabular}{c|rrr|rr|rr}
    \hline
          & \multicolumn{3}{c|}{Two-stage enumeration} & \multicolumn{2}{c|}{N-DDDR} & \multicolumn{2}{c}{N-DIDR} \\
          \hline
     Pattern & $(1,0,0)^{\mathsf T}$ & $(0,1,0)^{\mathsf T}$ & $(0,0,1)^{\mathsf T}$ & Obj.  & $\boldsymbol{x}$-sol.  & Obj.  & $\boldsymbol{x}$-sol. \\
     \hline
    1-1    & $\boldsymbol{-2160}$ & $-1800$  & $-1575$  & $-2160$ & $(1,0,0)^{\mathsf T}$ & \multirow{4}[0]{*}{$-1463$} & \multirow{4}[0]{*}{$(0,0,1)^{\mathsf T}$} \\
    1-2    & $\boldsymbol{-1800}$  & $\boldsymbol{-1800}$  & $\boldsymbol{-1800}$  & $-1800$ & $(0,0,1)^{\mathsf T}$ &       &  \\
    1-3     & $\boldsymbol{-1665}$ & $-1575$ & $-1485$ & $-1665$ & $(1,0,0)^{\mathsf T}$ &       &  \\
    1-4    & $-1800$  & $\boldsymbol{-2160}$ & $-1485$ & $-2160$ & $(0,1,0)^{\mathsf T}$ &       &  \\
    \hline
    \end{tabular}%
    }
  \label{tab:type1-2}%
\end{table}%

From Table \ref{tab:type1-2}, both the optimal solutions and objective values of the two-stage model by enumeration and N-DDDR are the same, confirming the finite convergence of the SDDiP algorithm. The model N-DDDR always yields a better objective value than the one of N-DIDR, indicating the benefits of considering decision-dependency. When we set $\lambda^S$-values the same, as shown in Pattern \#1-1, N-DDDR first builds the facility that has the highest impact on the mean values of demand, coinciding with our intuition that building such a facility will increase demand in later stages the most and as a result, it will bring the largest revenue. When we decrease all $\lambda^\mu$-values to $0.1$, N-DDDR chooses the facility with the highest $\lambda^S$-value, indicated in the optimal solution in Pattern \#1-3. In Patterns \#1-4, N-DDDR chooses the facility with $\lambda^{\mu}=0.9$ and $\lambda^S=0.5$. These results suggest that the impact on the first moment (e.g., mean values) plays a more important role than the impact on demand variance when choosing optimal facility-location solutions. 

\subsubsection{Results of Type 2 ambiguity set on two-stage instances}
Now we consider N-DDDR model with $T=2$ stages, $I=3$ facilities and $J=2$ customer sites. 
For Type 2 ambiguity set \eqref{eq:ambi_match}, assume that each facility has the same impact on different customer sites, i.e., $\lambda_{ji}^{\mu}=\lambda_i,\ \forall j\in[J]$. The empirical mean and covariance matrix are given by $\bar{\boldsymbol{\mu}}=(10,10)^{\mathsf T},\ \bar{\boldsymbol{\Sigma}} = \begin{pmatrix} 10&10\\10&10\end{pmatrix}$. Note that this type of ambiguity set is the most restricted one because it is defined by three equalities. As a result, we evaluate three different patterns with fixed $\lambda^\mu$- and $\sigma$-values given in Table \ref{tab:pattern_2}, which will make the ambiguity set \eqref{eq:ambi_match} non-empty. Table \ref{tab:type2-2} demonstrates the results of the two-stage model solved by enumeration, N-DDDR and N-DIDR solved by SDDiP under different patterns.
\begin{table}[ht!]
    \centering
        \caption{Patterns with varying $\lambda^\mu$-/$\sigma$-values that interpret how decisions affect mean/covariance.}
    \begin{tabular}{c|r|r}
    \hline
        Pattern & $(\lambda_{j1}^{\boldsymbol{\mu}},\lambda_{j2}^{\boldsymbol{\mu}},\lambda_{j3}^{\boldsymbol{\mu}})^{\mathsf T}$ & $(\sigma_1,\sigma_2,\sigma_3)^{\mathsf T}$ \\
        \hline
        2-1 & $(0.1,0.2,0.3)^{\mathsf T}$ & $(0.5,0.5,0.5)^{\mathsf T}$\\
        2-2 & $(0.1,0.2,0.3)^{\mathsf T}$ & $(0.9,0.5,0.1)^{\mathsf T}$\\
        2-3 & $(0.3,0.3,0.3)^{\mathsf T}$ & $(0.9,0.5,0.1)^{\mathsf T}$\\
        \hline
    \end{tabular}
    \label{tab:pattern_2}
\end{table}

\begin{table}[ht!]
  \centering
  \caption{Results of different models using Type 2 ambiguity sets}
   \resizebox{\textwidth}{!}{
    \begin{tabular}{c|rrr|rr|cc}
    \hline
          & \multicolumn{3}{c|}{Two-stage enumeration} & \multicolumn{2}{c|}{N-DDDR} & \multicolumn{2}{c}{N-DIDR} \\
          \hline
    Pattern  & $(1,0,0)^{\mathsf T}$ & $(0,1,0)^{\mathsf T}$ & $(0,0,1)^{\mathsf T}$ & Obj.  & $\boldsymbol{x}$-sol.  & \multicolumn{1}{r}{Obj.} & \multicolumn{1}{r}{$\boldsymbol{x}$-sol.} \\
    \hline
    2-1     & $-3780$ & $-3960$ & $\boldsymbol{-4140}$ & $-4140$ & $(0,0,1)^{\mathsf T}$ & \multirow{3}[0]{*}{$-3600$} & \multirow{3}[0]{*}{$(0,0,1)^{\mathsf T}$} \\
    2-2      & $-3780$ & $-3960$ & $\boldsymbol{-4140}$ & $-4140$ & $(0,0,1)^{\mathsf T}$ &       &  \\
    2-3    & $\boldsymbol{-4140}$ & $\boldsymbol{-4140}$ & $\boldsymbol{-4140}$ & $-4140$ & $(0,0,1)^{\mathsf T}$ &       &  \\
    \hline
    \end{tabular}%
    }
  \label{tab:type2-2}%
\end{table}%

From Table \ref{tab:type2-2}, in Pattern \#2-1, when all the $\sigma$-values are the same, N-DDDR builds the facility with the highest impact on the mean. In Pattern \#2-2, when the third facility has the highest impact on the mean ($\lambda_{j3}^{\mu}=0.3$) and the lowest impact on the covariance matrix ($\sigma_3=0.1$), N-DDDR still builds the third one, indicating the importance of mean values of demand. 

\subsubsection{Results of Type 3 ambiguity set on two-stage instances}
For Type 3 ambiguity set \eqref{eq:ambi_bound}, we set bounding parameters as $\gamma=1000,\ \eta=500$, the empirical mean and covariance matrix as $\bar{\boldsymbol{\mu}}=(10,10)^{\mathsf T},\ \bar{\boldsymbol{\Sigma}}=\begin{pmatrix}0.1 & 0.2\\0.2 &0.9\end{pmatrix}$. 
We evaluate four different patterns with fixed $\lambda^\mu$- and $\sigma$-values given in Table \ref{tab:pattern_3}. Then in Table \ref{tab:type3-3}, we show the results of the two-stage model solved by enumeration, N-DDDR and N-DIDR for different patterns given in Table \ref{tab:pattern_3}. 
\begin{table}[ht!]
    \centering
        \caption{Patterns with varying $\lambda^\mu$-/$\sigma$-values that interpret how decisions affect mean/covariance.}
    \begin{tabular}{c|r|r}
    \hline
        Pattern & $\begin{pmatrix}\lambda_{11}^{\boldsymbol{\mu}},\lambda_{12}^{\boldsymbol{\mu}},\lambda_{13}^{\boldsymbol{\mu}}\\
        \lambda_{21}^{\boldsymbol{\mu}},\lambda_{22}^{\boldsymbol{\mu}},\lambda_{23}^{\boldsymbol{\mu}}\end{pmatrix}$ & $(\sigma_{1},\sigma_2,\sigma_3)^{\mathsf T}$ \\
        \hline
        3-1 &  $\begin{pmatrix}0.1,0.5,0.9\\0.1,0.5,0.9 \end{pmatrix}$ & $(0.5,0.5,0.5)^{\mathsf T}$\\
        3-2 &  $\begin{pmatrix}0.5,0.5,0.5\\0.5,0.5,0.5 \end{pmatrix}$& $(0.9,0.5,0.1)^{\mathsf T}$\\
        3-3 &  $\begin{pmatrix}0.1,0.5,0.9\\0.1,0.5,0.9 \end{pmatrix}$ & $(0.1,0.5,0.9)^{\mathsf T}$\\
        3-4 & $\begin{pmatrix}0.1,0.5,0.9\\0.9,0.5,0.1 \end{pmatrix}$ & $(0.5,0.5,0.5)^{\mathsf T}$\\
        \hline
    \end{tabular}
    \label{tab:pattern_3}
\end{table}

\begin{table}[ht!]
  \centering
  \caption{Results of different models under Type 3 ambiguity set and correlated demand}
   \resizebox{\textwidth}{!}{
    \begin{tabular}{c|rrr|rr|cc}
    \hline
          & \multicolumn{3}{c|}{Two-stage exact} & \multicolumn{2}{c|}{N-DDDR} & \multicolumn{2}{c}{N-DIDR} \\
          \hline
    Pattern  & $(1,0,0)^{\mathsf T}$ & $(0,1,0)^{\mathsf T}$ & $(0,0,1)^{\mathsf T}$ & Obj.  & $\boldsymbol{x}$-sol.  & \multicolumn{1}{r}{Obj.} & \multicolumn{1}{r}{$\boldsymbol{x}$-sol.} \\
    \hline
    3-1     & $-2700$ & $-3150$ & $\boldsymbol{-3856.2}$ & $-3856.2$ & $(0,0,1)^{\mathsf T}$ & \multirow{4}[0]{*}{$-2701$} & \multirow{4}[0]{*}{$(1,0,0)^{\mathsf T}$} \\
    3-2      & $-2950$ & $-3150$ & $\boldsymbol{-3350}$ & $-3350$ & $(0,0,1)^{\mathsf T}$ &       &  \\
    3-3    & $-2700$ & $-3150$ & $\boldsymbol{-3625.4}$ & $-3625.4$ & $(0,0,1)^{\mathsf T}$ &       &  \\
    3-4    & $-2700$ & $-3150$ & $\boldsymbol{-4201.8}$ & $-4201.8$ & $(0,0,1)^{\mathsf T}$ &       &  \\
    \hline
    \end{tabular}%
    }
  \label{tab:type3-3}%
\end{table}%

From Table \ref{tab:type3-3}, when the $\lambda^{\mu}$-values are the same as shown in Pattern \#3-2, N-DDDR builds the facility with the lowest impact on the covariance matrix, which is different from the previous two ambiguity sets. In other  patterns, N-DDDR always builds the facility with the highest impact on the mean values of demand for both customer locations. When different facilities have the highest impact on the demand in the two locations, the location with smaller demand variance will play a more important role in choosing facilities to build.

{
\subsection{Results of Larger Instances under Type 1 Ambiguity Set}
\label{sec:compu_type1}
We first consider N-DDDR model with $T=3$ stages, $I=10$ facilities, $J=20$ customer sites and Type 1 ambiguity set \eqref{eq:ambi_ineq} in Section \ref{sec:ambiguity_1}. We set the bounding parameters $\epsilon^{\mu}=25,\ \underline{\epsilon}^S=0.1,\ \bar{\epsilon}^S=1.9$.  Parameters $\lambda_{ji}^{\mu},\ \lambda_{ji}^S$ follow exponential functions in terms of the distance between customer site $j$ and facility $i$ so that farther facilities have lower impacts on the first and second moments of the demand, i.e., $\lambda_{ji}^{\mu}= e^{-\textrm{dist}(i,j)/25},\ \lambda_{ji}^S= e^{-\textrm{dist}(i,j)/50}$ for all $i\in[I],\ j\in[J]$, and then they are normalized to ensure that the sum of impacts over all facilities equals to 1, i.e., $\sum_{i\in[I]}\lambda_{ji}^{\mu}=\sum_{i\in[I]}\lambda_{ji}^S=1,\ \forall j\in[J]$. 

We sample $K$ data points following $\mathcal{N}(\bar{\mu}_j, \bar{\sigma}_j^2)$ to construct the discrete support $\Xi_t$ for each $t\in[T]$, and set the demand variation coefficient $\bar{\rho}=\bar{\sigma}_j/\bar{\mu}_j$ to 0.8 for each $j\in[J]$ by default, where we vary it in Section \ref{sec:variance}.  We then apply SDDiP algorithm to solve both the N-DDDR and N-DIDR with all $\lambda^\mu$- and $\lambda^S$-values set to 0. The locations of potential facilities and customer sites are displayed in Figure \ref{fig:locations}, where triangles represent customer sites and circles stand for potential facilities.
\begin{figure}[ht!]
\centering
\begin{tikzpicture}[scale=0.8]
\begin{axis}[legend style={nodes={scale=0.8, transform shape}}, 
legend pos=north east]
    \addplot[
        scatter/classes={       a={mark=triangle*, fill=red},
        b={mark=*,fill=blue}},
        scatter, mark=*, only marks, 
        scatter src=explicit symbolic,
        nodes near coords*={\Label},
        visualization depends on={value \thisrow{label} \as \Label} 
    ] table [meta=class] {
        x y class label
        81	99 a {}
        96	75 a {}
        31	28 a  {}
        68	77 a {}
        85	10 a {}
        99	43 a {}
        8	86 a {}
        4	98 a {}
        16	27 a {}
        80	12 a {}
        9	2 a {}
        38	61 a {}
        95	19 a {}
        47	24 a {}
        60	95 a {}
        27	5 a {}
        58	49 a {}
        70	13 a {}
        2	84 a {}
        61	57 a {}
        42	42 b \#1
        72	68 b \#2
        1	21 b \#3
        30	86 b \#4
        15	3  b \#5
        9	64 b \#6
        18	40 b \#7 
        33	52 b \#8
        37	13 b \#9
        50	19 b \#10
    };
    \legend{Customer sites, Potential facilities}
\end{axis}
\end{tikzpicture}
\caption{Locations of customer sites and potential facilities on a 100$\times$100 grid} \label{fig:locations}
\end{figure}

\subsubsection{Objective values with different support sizes}
We vary the number $K$ of data samples in the discrete, finite support from 10 to 100 and display the objective values of N-DDDR and N-DIDR in Figure \ref{fig:trend}(a), respectively, where Figure \ref{fig:trend}(b) zooms in Figure \ref{fig:trend}(a) by dropping the unbounded cases. 
\begin{figure}[ht!]
    \centering
    \subfigure[The original result]{
    \resizebox{0.45\textwidth}{!}{%
\begin{tikzpicture}
  \begin{axis}
  [
    xlabel={Support size $K$},
    ylabel={Objective values},
    xtick={10,20,30,40,50,60,70,80,90,100},
    legend pos=south east,
]
    \addplot coordinates {
      (10,  -1050940)
      (20,  -1050940)
      (30,  -165161.6351)
      (40,  -164660.2631)
      (50,  -164600.1601)
      (60,  -164237.1834)
      (70,  -164237.1807)
      (80,  -164237.1799)
      (90,  -163907.1679)
      (100, -163907.4558)
    };
    \addplot coordinates {
      (10,	-150155.3816)
      (20,  -147950.8925)
      (30,  -147949.4927)
      (40,  -147950.8925)
      (50,  -147949.5705)
      (60,  -147120.1563)
      (70,  -147120.1563)
      (80,  -147120.3647)
      (90,  -146534.9203)
      (100, -146534.7342)
    };
    \legend{N-DDDR, N-DIDR}
  \end{axis}
\end{tikzpicture}%
}
}
\subfigure[A truncated version by dropping all the  unbounded cases]{
\resizebox{0.45\textwidth}{!}{%
\begin{tikzpicture}
  \begin{axis}
  [
    xlabel={Support size $K$},
    ylabel={Objective values},
    xtick={10,20,30,40,50,60,70,80,90,100},
    legend pos=south east,
]
    \addplot coordinates {
      (30,  -165161.6351)
      (40,  -164660.2631)
      (50,  -164600.1601)
      (60,  -164237.1834)
      (70,  -164237.1807)
      (80,  -164237.1799)
      (90,  -163907.1679)
      (100, -163907.4558)
    };
    \addplot coordinates {
      (30,  -147949.4927)
      (40,  -147950.8925)
      (50,  -147949.5705)
      (60,  -147120.1563)
      (70,  -147120.1563)
      (80,  -147120.3647)
      (90,  -146534.9203)
      (100, -146534.7342)
    };
  \end{axis}
\end{tikzpicture}%
}
}
    \caption{Objective values of N-DDDR and N-DIDR under Type 1 ambiguity set and varying support sizes $K$}
    \label{fig:trend}
\end{figure}

From Figure \ref{fig:trend}(a), when $K=10,\ 20$, the N-DDDR model is unbounded with an empty ambiguity set \eqref{eq:ambi_ineq}, mainly due to a lack of data points in the discrete support. By increasing the support size $K$, the objective values of N-DDDR increase. Recall that the worst-case scenario is calculated by the inner maximization problem, and therefore, larger-sized discrete supports lead to higher worst-case objectives. Overall, we are minimizing the N-DDDR objective function, and thus lower objective values are more favorable. More data points $K$ in the discrete support can either be interpreted as a more risk-averse altitude, or represent a better approximation of the continuous distribution. From Figure \ref{fig:trend}(b), the objective values of N-DDDR and N-DIDR both have step-wise increments.  That is, when we include more data points, the objective values may stay constant or take a step upward, depending on whether the inclusion of these data points changes the worst-case scenarios. Moreover, N-DDDR always yields better objective values than N-DIDR, indicating the benefits of considering decision-dependency. 

\subsubsection{Objective values with different sample variance and distributions}\label{sec:variance}
Next, we fix the support size $K=100$ and vary the demand variation coefficient $\bar{\rho}=\bar{\sigma}_j/\bar{\mu}_j$ from 0.2 to 1 for all $j\in[J]$. To further illustrate the impact of demand variations on the objective values, we also compare the results of different distributions of which the data points come from. Figure \ref{fig:variation} displays the objective values of N-DDDR with varying demand variation coefficients $\bar{\rho}$ and Normal/Log-normal distributions, respectively, where we drop the demand variations that make the problem unbounded (i.e., make the ambiguity sets empty). To be comparable with Normal distributions, we set the scale parameter (the median of the Log-normal distribution) to be the empirical mean of the Normal distribution, i.e., $\bar{\mu}_j$, the location parameter (parameter $\mu$ of the Log-normal distribution) to be $\log(\bar{\mu}_j)$ and the shape parameter (parameter $\sigma$ of the Log-normal distribution) to be $\bar{\rho}\log(\bar{\mu}_j)$ for each $j\in[J]$. 

\begin{figure}[ht!]
    \centering
    \subfigure[Normal distribution]{
    \resizebox{0.45\textwidth}{!}{%
\begin{tikzpicture}
  \begin{axis}
  [
    xlabel={Demand variation $\bar{\rho}$},
    ylabel={Objective values},
    xtick={0.2, 0.3, 0.4, 0.5, 0.6, 0.7, 0.8, 0.9, 1},
    legend style={at={(0.6,0.5)},anchor=west}
]
    \addplot coordinates {
    (0.6,	-161526.6617)
    (0.7,	-162497.953)
    (0.8,	-163907.4558)
    (0.9,	-165677.674)
    (1,	-167857.0316)
    };
  \end{axis}
\end{tikzpicture}%
}
}
\subfigure[Log-normal distribution]{
\resizebox{0.45\textwidth}{!}{%
\begin{tikzpicture}
  \begin{axis}
  [
    xlabel={Demand variation $\bar{\rho}$},
    ylabel={Objective values},
    xtick={0.2, 0.3, 0.4, 0.5, 0.6, 0.7, 0.8, 0.9, 1},
    legend pos=north west
]
    \addplot coordinates {
    (0.2,	-163861.1416)
    (0.3,	-164544.1803)
    (0.4,	-168819.5919)
    };
  \end{axis}
\end{tikzpicture}%
}
}
    \caption{Objective values of N-DDDR under Type 1 ambiguity set, varying demand variations $\bar{\rho}$ and distributions}
    \label{fig:variation}
\end{figure}

In Figure \ref{fig:variation}, the objective values with Normal and Log-normal distributions have totally different behaviors with respect to demand variations. When $\bar{\rho}$ is low (i.e., $\bar{\rho}\le 0.5$), the ambiguity sets with discrete supports constructing by Normal distributions are empty, because the data points in the discrete support mostly concentrate around the empirical mean and lack of diversity. On the contrary, the problem with Log-normal distributions becomes unbounded when $\bar{\rho}$ is high (i.e., $\bar{\rho}\ge 0.5$). This is because of the long-tail characteristic of Log-normal distributions. Under increasing demand variations, it is more likely to include extreme scenarios in the discrete support when sampling from a Log-normal distribution, and having too many deviated data points from the empirical mean is hard to construct a non-empty ambiguity set \eqref{eq:ambi_ineq}. It is also worth noting that the objective values with Normal and Log-normal distributions both decrease as demand variation increases.

\subsubsection{Optimal solutions with varying budgets and transportation costs}
We fix the support size $K=100$, demand variation $\bar{\rho}$ at 0.8, and increase the building budget $N$ from 100 to 500. Table \ref{tab:budget1} displays the optimal objective values and solutions of models N-DDDR and N-DIDR with varying budgets, respectively.
\begin{table}[ht!]
  \centering
  \caption{Optimal solutions of N-DDDR and N-DIDR with varying budgets}
    \begin{tabular}{rrrrr}
    \hline
    Budget $N$ & N-DDDR Obj.   & N-DDDR Sol. & N-DIDR Obj. & N-DIDR Sol. \\
    \hline
    100     & $-163,907$ & $1$ & $-146,535$ & $1$ \\
    300     & $-169,148$ & $[2,4,10]$ & $-150,479$ & $[2,6,10]$\\
    500     & $-171,986$ & $[2,4,6,8,10]$ & $-151,323$ & $[1,2,5,6,10]$ \\
    \hline
    \end{tabular}%
  \label{tab:budget1}%
\end{table}%

In Table \ref{tab:budget1}, when we only have budgets to build one facility at the first stage, the optimal solutions of N-DDDR and N-DIDR both choose facility \#1. Combining with Figure \ref{fig:locations}, facility \#1 is in the most central location. With higher budget values $N=300$ and $N=500$, the optimal solutions of N-DDDR do not include facility \#1 anymore and the objective values get improved by building more facilities. Moreover, N-DDDR always yields better objective values than N-DIDR by building facilities having bigger impacts on the demand mean.

To not take relative locations into account, we set all the transportation costs to 10, and record the optimal objective values and solutions in Table \ref{tab:budget2}. We also display the impacts on the first and second moments of all customer sites by calculating $\sum_{j\in [J]}\lambda_{ji}^{\mu},\ \sum_{j\in [J]}\lambda_{ji}^{S}$ for each $i\in[I]$ in Table \ref{tab:impact}.

\begin{table}[ht!]
  \centering
  \caption{Optimal solutions of N-DDDR and N-DIDR with varying budgets and fixed transportation cost}
    \begin{tabular}{rrrrr}
    \hline
    Budget $N$ & N-DDDR Obj.   & N-DDDR Sol. & N-DIDR Obj. & N-DIDR Sol. \\
    \hline
    100     & $-162,752$ & $2$ & $-145,811$ & $9$\\
    300     & $-165,461$ & $[2,4,8]$ & $-145,815$ & $[7,8,9]$\\
    500     & $-168,114$ & $[1,2,4,8,10]$ & $-145,820$ & $[4,5,7,8,10]$\\
    \hline
    \end{tabular}%
  \label{tab:budget2}%
\end{table}%

\begin{table}[htbp]
  \centering
  \caption{Total impact on the first and second moments of each facility}
    \begin{tabular}{rrrrrrrrrrr}
    \hline
    $I$     & \#1     & \#2     & \#3     & \#4     & \#5     & \#6     & \#7     & \#8     & \#9     & \#10 \\
    \hline
    $\sum_{j\in [J]}\lambda_{ji}^{\mu}$ & 2.11  & 3.47  & 1.08  & 2.15  & 1.43  & 1.69  & 1.47  & 1.86  & 2.09  & 2.65 \\
    $\sum_{j\in [J]}\lambda_{ji}^{S}$ & 2.23  & 2.58  & 1.48  & 2.02  & 1.57  & 1.80  & 1.85  & 2.11  & 2.06  & 2.31 \\
    \hline
    \end{tabular}%
  \label{tab:impact}%
\end{table}%

From Tables \ref{tab:budget2} and \ref{tab:impact}, facility \#2 has the largest total impact on the first and second moments of the uncertain demand among all facilities, which is built in the optimal solutions of N-DDDR. However, without the decision-dependency settings, the optimal solutions of N-DIDR always choose facilities having smaller impacts on the first and second moments, leading to worse objective values than N-DDDR.

\subsubsection{Computational time}
We first compare the computational time of models N-DDDR and N-DIDR under Type 1  ambiguity set. We first fix $I=10,\ J=20,\ T=3$ and vary the number $K$ of data points in the support from 10 to 100. The computational time results are displayed in Figure \ref{fig:time1}(a). Then we fix $K=10,\ T=3$ and vary the number $I$ of facilities from 10 to 50 while setting $J=2I$ in Figure \ref{fig:time1}(b). Finally, we fix $K=10, \ I=10,\ J=20$ and vary the number $T$ of stages from 3 to 8 in Figure \ref{fig:time1}(c).
\begin{figure}[ht!]
    \centering
    \subfigure[different support sizes $K$]{
 \resizebox{0.45\textwidth}{!}{%
\begin{tikzpicture}
  \begin{axis}
  [
    xlabel={Support size $K$},
    ylabel={Computational time (sec.)},
    xtick={10,20,30,40,50,60,70,80,90,100},
    legend pos=south east,
]
    \addplot coordinates {
      (10, 31.58350121)
      (20, 24.52017935)
      (30, 76.64471387)
      (40, 95.17694093)
      (50, 131.9146672)
      (60, 130.9996011)
      (70, 167.3638442)
      (80, 212.6122302)
      (90, 335.08)
      (100, 270.6766078)
    };
    \addplot coordinates {
      (10,	14.21529542)
      (20,  28.26336449)
      (30,  58.89553567)
      (40,  72.37103063)
      (50,  99.56892536)
      (60,  104.3380809)
      (70,  133.3971975)
      (80,  161.5314503)
      (90,  201.8281727)
      (100, 191.2663862)
    };
    \legend{N-DDDR, N-DIDR}
  \end{axis}
\end{tikzpicture}%
}
}
\subfigure[different numbers of facilities $I$]{
 \resizebox{0.45\textwidth}{!}{%
\begin{tikzpicture}
  \begin{axis}
  [
    xlabel={Number of facilities $I$},
    ylabel={Computational time (sec.)},
    xtick={10,20,30,40,50},
    legend pos=south east,
]
    \addplot coordinates {
      (10, 31.58350121)
      (20, 130.3365129)
      (30, 736.351985)
      (40, 1326.05987)
      (50, 2626.964099)
    };
    \addplot coordinates {
      (10,	14.21529542)
      (20, 73.16349176)
      (30, 267.6631899)
      (40, 538.7447037)
      (50, 1034.840937)
    };
    \legend{N-DDDR, N-DIDR}
  \end{axis}
\end{tikzpicture}%
}
}
\subfigure[different stages $T$]{
 \resizebox{0.45\textwidth}{!}{%
\begin{tikzpicture}
  \begin{axis}
  [
    xlabel={Time stages $T$},
    ylabel={Computational time (sec.)},
    xtick={3,4,5,6,7,8},
    legend pos=south east,
]
    \addplot coordinates {
      (3, 120.1155771)
      (4, 260.2855509)
      (5, 363.7509949)
      (6, 479.5241542)
      (7, 610.2713716)
      (8, 693.8911787)
    };
    \addplot coordinates {
      (3,	56.73110991)
      (4, 173.8702674)
      (5, 238.3046549)
      (6, 336.0905426)
      (7, 418.439642)
      (8, 479.7497408)
    };
    \legend{N-DDDR, N-DIDR}
  \end{axis}
\end{tikzpicture}%
}
}
    \caption{Computational time of N-DDDR and N-DIDR under Type 1  ambiguity set and different support sizes $K$, numbers of facilities $I$, and stages $T$}
    \label{fig:time1}
\end{figure}

In Figure \ref{fig:time1}, the computational time increases approximately linearly with respect to the support size $K$ and the number $T$ of stages, while it increases exponentially with respect to the numbers of facilities $I$ and customer sites $J$, due to the existence of McCormick constraints. Moreover, the N-DDDR model is always more time-consuming than the N-DIDR counterpart, although it has superior performance in terms of objective values as we note before. 

\subsection{Results of Larger Instances under Type 3 Ambiguity Set}\label{sec:compu_type3}
The default setting of the N-DDDR model in this section has $T=3$ stages, $I=3$ facilities and $J=6$ customer sites. For Type 3 ambiguity set \eqref{eq:ambi_bound}, we set bounding parameters $\gamma=10,\ \eta=100$. Parameters $\lambda_{ji}^{\mu}$ are the same as described in Section \ref{sec:compu_type1}, and ${\lambda}_{i}^{cov}$ are drawn uniformly between 0 and 1 for all $i\in[I],\ j\in[J]$. We then normalize parameters $\lambda_{ji}^{\mu},\ {\lambda}_{i}^{cov}$ to ensure that the sum over all facilities equals to 1, i.e., $\sum_{i\in[I]}\lambda_{ji}^{\mu}=\sum_{i\in[I]}{\lambda}_{i}^{cov}=1,\ \forall j\in[J]$. 

We sample $K$ data points following $\mathcal{N}(\bar{\mu}_j,{\bar{\sigma}_j^2})$ for each $j\in[J]$ to construct the discrete support, and set the demand variation coefficient $\bar{\rho}=\bar{\sigma}_j/\bar{\mu}_j$ to 0.8 for all $j\in[J]$ by default, where we vary it in Section \ref{sec:variance2}. Then the empirical covariance $\bar{\boldsymbol{\Sigma}}$ is set to the sample covariance matrix of the $K\times J$ data points among all customer sites. The locations of potential facilities and customer sites are displayed in Figure \ref{fig:locations2}, where triangles represent customer sites and circles represent potential facilities.

\begin{figure}[ht!]
\centering
\begin{tikzpicture}[scale=0.8]
\begin{axis}[legend style={nodes={scale=0.8, transform shape}}, 
legend pos=north east]
    \addplot[
        scatter/classes={       a={mark=triangle*, fill=red},
        b={mark=*,fill=blue}},
        scatter, mark=*, only marks, 
        scatter src=explicit symbolic,
        nodes near coords*={\Label},
        visualization depends on={value \thisrow{label} \as \Label} 
    ] table [meta=class] {
        x y class label
        19	21 a {}
        35	87 a {}
        39	3 a {}
        53	66 a {}
        41	41 a {}
        66	54 a {}
         42	31 b \#1
        72	15 b \#2
        1	10 b \#3
    };
    \legend{Customer sites, Potential facilities}
\end{axis}
\end{tikzpicture}
\caption{Locations of customer sites and potential facilities on a 100$\times$100 grid} 
\label{fig:locations2}
\end{figure}

\subsubsection{Objective values with different support sizes}
We vary the values of $K$ from 10 to 100 and display bounds on the objective values of models N-DDDR and N-DIDR in Figures \ref{fig:trend2}, respectively, where ``LB'' indicates valid lower bounds using Relaxed Lagrangian Cuts introduced in Section \ref{sec:relax}, and ``UB'' stands for valid upper bounds provided by the inner approximation scheme in Section \ref{sec:upper} with $\boldsymbol{U},\boldsymbol{V}$ being identity matrices. 

\begin{figure}[ht!]
    \centering
    \subfigure[Bounds on objective values of N-DDDR]{
    \resizebox{0.45\textwidth}{!}{%
\begin{tikzpicture}
  \begin{axis}
  [
    xlabel={Support size $K$},
    ylabel={Objective values},
    xtick={10,20,30,40,50},
    legend pos=south east,
]
    \addplot coordinates {
    (10,-36218.20355)
    (20, -35703)
    (30, -34756.28405)
    (40, -34951.2362) 
    (50, -35678.14139)
    };
    \addplot coordinates {
    (10, -35692.00035)
    (20, -35693.00009)
    (30, -34415)
    (40, -34415)
    (50, -34414)
    };
    \legend{N-DDDR-LB, N-DDDR-UB}
  \end{axis}
\end{tikzpicture}%
}
}
\subfigure[Bounds on objective values of N-DIDR]{
    \resizebox{0.45\textwidth}{!}{%
\begin{tikzpicture}
  \begin{axis}
  [
    xlabel={Support size $K$},
    ylabel={Objective values},
    xtick={10,20,30,40,50},
    legend pos=south east,
]
    \addplot coordinates {
    (10, -35703.0001)
    (20, -35703.00001)
    (30, -34417.08333)
    (40, -34419.00009)
    (50, -34419.00008)
    };
        \addplot coordinates {
    (10, -35693)
    (20, -35683.64286)
    (30, -34413)
    (40, -34411.58333)
    (50, -34411.58333)
    };
    \legend{N-DIDR-LB, N-DIDR-UB}
  \end{axis}
\end{tikzpicture}%
}
}
    \caption{Objective values of N-DDDR and N-DIDR under Type 3 ambiguity set and varying support sizes $K$}
    \label{fig:trend2}
\end{figure}

In Figure \ref{fig:trend2}, the objective values of N-DDDR's UB and N-DIDR's LB and UB all increase stepwise with increased support sizes $K$, and the objective values of N-DIDR are slightly higher than N-DDDR's UB. It is also worth noting that the relative gaps of N-DDDR are always within 4\% while the scale of the relative gaps of N-DIDR is at $10^{-4}$, showing the close proximity of LB and UB provided by our algorithms. Moreover, both the LB and UB of N-DDDR and N-DIDR choose to build facility \#1 in the first stage of the optimal solutions, which locates centrally and also has the largest impact on the mean and covariance of the uncertain demand.

\subsubsection{Objective values with different sample variance and distributions}\label{sec:variance2}
Next we fix $K=10$ data points in the support and vary the demand variation $\bar{\rho}=\bar{\sigma}_j/\bar{\mu}_j$ from 0.2 to 1. Figure \ref{fig:variation2} displays the objective values of model N-DDDR's LB and UB with respect to Normal and Log-normal distributions, respectively. We only display the demand variations that make the ambiguity sets non-empty and drop the unbounded cases. 
\begin{figure}[ht!]
    \centering
    \subfigure[Normal distributions]{
    \resizebox{0.45\textwidth}{!}{%
\begin{tikzpicture}
  \begin{axis}
  [
    xlabel={Demand variation $\bar{\rho}$},
    ylabel={Objective values},
    xtick={0.5,0.6,0.7,0.8,0.9,1},
    legend pos=south east,
]
    \addplot coordinates {
    (0.5, -42640.08482)
    (0.6, -39838.42691)
    (0.7, -37793.91049)
    (0.8, -36218.20355)
    (0.9, -36357.51677)
    (1, -36964.35362)
    };
    \addplot coordinates {
    (0.5, -35593)
    (0.6, -35382)
    (0.7, -35533)
    (0.8, -35692.00035)
    (0.9, -36214)
    (1, -36912)
    };
    \legend{N-DDDR-LB, N-DDDR-UB}
  \end{axis}
\end{tikzpicture}%
}
}
\subfigure[Log-normal distributions]{
\resizebox{0.45\textwidth}{!}{%
\begin{tikzpicture}
  \begin{axis}
  [
    xlabel={Demand variation $\bar{\rho}$},
    ylabel={Objective values},
    xtick={0.2,0.3,0.4,0.5,0.6,0.7,0.8},
    legend pos=south east,
]
    \addplot coordinates {
    (0.2, -35941)
    (0.3, -34129)
    (0.4, -31575)
    (0.5, -30497)
    (0.6, -29215)
    (0.7, -28489)
    };
    \addplot coordinates {
    (0.2, -35927)
    (0.3, -34126)
    (0.4, -31574)
    (0.5, -30496)
    (0.6, -29215)
    (0.7, -28489)
    };
    \legend{N-DDDR-LB, N-DDDR-UB}
  \end{axis}
\end{tikzpicture}%
}
}
    \caption{Objective values of N-DDDR under Type 3  ambiguity set and varying demand variations $\bar{\rho}$}
    \label{fig:variation2}
\end{figure}

In Figure \ref{fig:variation2}, similarly, when the demand variation $\bar{\rho}$ is low (i.e., $\bar{\rho}\le 0.4$), the ambiguity sets constructing by Normal distributions become empty, while the ones constructing by Log-normal distributions become empty when the demand variation is high (i.e., $\bar{\rho}\ge 0.8$). Moreover, the gaps between LB and UB decrease as demand variation increases with Normal distributions, while the gaps are significantly reduced with Log-normal distributions. 

\subsubsection{Computational time}\label{sec:time}
Lastly, we compare the computational time of solving models N-DDDR and N-DIDR under Type 3  ambiguity set.  We first fix $I=3,\ J=6,\ T=3$ and vary the support size $K$ from 10 to 50, displayed in Figure \ref{fig:time2}(a). Then we fix $K=10,\ T=3$ and vary the number $I$ of facilities from 3 to 6 while setting $J=2I$ in Figure \ref{fig:time2}(b). Finally, we fix $K=10, \ I=3,\ J=6$ and vary the number $T$ of stages from 3 to 6 in Figure \ref{fig:time2}(c). The time limit for solving each instance is set as 7200 seconds.
\begin{figure}[ht!]
    \centering
    \subfigure[different support sizes $K$]{
 \resizebox{0.45\textwidth}{!}{%
\begin{tikzpicture}
  \begin{axis}
  [
    xlabel={Support size $K$},
    ylabel={Computational time (sec.)},
    xtick={10,20,30,40,50},
    legend pos=north west,
]
    \addplot coordinates {
      (10, 1774.798429)
      (20, 2543.330376)
      (30, 3227.848633)
      (40, 4521.001386)
      (50, 5851.10631)
    };
    \addplot coordinates {
      (10,	105.89)
      (20, 174.8605071)
      (30, 323.0178448)
      (40, 348.3627082)
      (50, 440.8941127)
    };
        \addplot coordinates {
      (10,	189.38)
      (20, 558.3)
      (30, 1069.837829)
      (40, 2092.6)
      (50, 3810.117868)
    };
        \addplot coordinates {
      (10,	105.1502885)
      (20, 442.0286124)
      (30, 652.365994)
      (40, 1845.725809)
      (50, 2979.171281)
    };
    \legend{N-DDDR-LB, N-DIDR-LB,N-DDDR-UB, N-DIDR-UB}
  \end{axis}
\end{tikzpicture}%
}
}
\subfigure[different numbers of facilities $I$]{
 \resizebox{0.45\textwidth}{!}{%
\begin{tikzpicture}
  \begin{axis}
  [
    xlabel={Number of facilities $I$},
    ylabel={Computational time (sec.)},
    xtick={3,4,5,6},
    legend pos=south east,
]
    \addplot coordinates {
      (3, 1732.4)
      (4, 7310.895272)
    };
    \addplot coordinates {
      (3,105.89)
       (4, 4123.2)
    };
        \addplot coordinates {
      (3,189.38)
      (4, 951.35)
      (5, 3311.3)
    };
        \addplot coordinates {
      (3,105.1502885)
      (4, 474.85)
      (5, 1153)
      (6, 2813.5)
    };
  \end{axis}
\end{tikzpicture}%
}
}
\subfigure[different stages $T$]{
 \resizebox{0.45\textwidth}{!}{%
\begin{tikzpicture}
  \begin{axis}
  [
    xlabel={Time stages $T$},
    ylabel={Computational time (sec.)},
    xtick={3,4,5,6},
    legend pos=south east,
]
    \addplot coordinates {
      (3, 1616.3)
      (4, 3580)
      (5, 5118.2)
    };
    \addplot coordinates {
      (3, 124.7516466)
      (4, 212.0832087)
      (5, 307.6340695)
      (6, 388.8025014)
    };
        \addplot coordinates {
      (3, 189.38)
      (4, 227.46)
      (5, 299.34)
      (6, 381.29)
    };
        \addplot coordinates {
      (3, 36.66909947) 
      (4, 61.9611286)
      (5, 118.0180072)
      (6, 213.2883482)
    };
  \end{axis}
\end{tikzpicture}%
}
}
    \caption{Computational time of N-DDDR and N-DIDR under Type 3  ambiguity set and different support sizes $K$, number of facilities $I$, and stages $T$}
    \label{fig:time2}
\end{figure}

In Figure \ref{fig:time2}, Type 3 ambiguity set \eqref{eq:ambi_bound} makes N-DDDR more difficult to solve than Type 1 ambiguity set. Comparing different approximation schemes for solving model N-DDDR under Type 3 ambiguity set, UB is the fastest as it solves a stage-wise MILP in both forward and backward step, LB is the most time-consuming as it solves a stage-wise MISDP in both forward and backward steps.
}

\section{Conclusions}
\label{sec:con}
In this paper, we studied multistage mixed-integer DRO model with decision-dependent moment-based ambiguity sets. We also extended the models to risk-averse cases by replacing the expectation with a coherent risk measure in the objective function. We recast the two problems as multistage stochastic MILP/MISDP and applied variants of SDDiP to solve them. Via numerical studies, we showed that N-DDDR always yielded a better objective value than that of its decision-independent counterpart. {Also, our solution approaches converged to the true optimal results under Types 1 and 2 ambiguity sets, and yielded small gaps between lower- and upper-bounds for N-DDDR under Type 3 ambiguity set.} 

The ambiguity sets used in this paper are all moment-based.  However, this  ambiguity sets do not have asymptotic consistency, i.e., we can not recover the true optimal objective value of the stochastic program as the number of data points increases to infinity. Therefore, it will be interesting to construct ambiguity sets based on some divergence measures, such as Wasserstein metric, and extend  such sets for the decision-dependent setting in our future research studies. 

~\\
{\bf Acknowledgements:}
The authors sincerely thank the Associate Editor and two reviewers for their helpful review and constructive feedback. 
The authors are grateful for the support from the United States National Science Foundation Grants \#1727618, \#1709094, and Department of Engineering (DoE) grant \#DE-SC0018018 for this project.

%
%



\begin{thebibliography}{44}
\providecommand{\natexlab}[1]{#1}
\providecommand{\url}[1]{{#1}}
\providecommand{\urlprefix}{URL }
\expandafter\ifx\csname urlstyle\endcsname\relax
  \providecommand{\doi}[1]{DOI~\discretionary{}{}{}#1}\else
  \providecommand{\doi}{DOI~\discretionary{}{}{}\begingroup
  \urlstyle{rm}\Url}\fi
\providecommand{\eprint}[2][]{\url{#2}}

\bibitem[{Ahmadi and Hall(2017)}]{ahmadi2017sum}
Ahmadi AA, Hall G (2017) Sum of squares basis pursuit with linear and second
  order cone programming. Algebraic and {G}eometric {M}ethods in {D}iscrete
  {M}athematics 685:27--53

\bibitem[{Basciftci et~al.(2019)Basciftci, Ahmed, and
  Shen}]{basciftci2019distributionally}
Basciftci B, Ahmed S, Shen S (2019) Distributionally robust facility location
  problem under decision-dependent stochastic demand. arXiv preprint
  arXiv:191205577

\bibitem[{Ben-Tal et~al.(2009)Ben-Tal, El~Ghaoui, and
  Nemirovski}]{ben2009robust}
Ben-Tal A, El~Ghaoui L, Nemirovski A (2009) Robust Optimization, vol~28.
  Princeton University Press

\bibitem[{Bertsimas et~al.(2011)Bertsimas, Brown, and
  Caramanis}]{bertsimas2011theory}
Bertsimas D, Brown DB, Caramanis C (2011) Theory and applications of robust
  optimization. SIAM Review 53(3):464--501

\bibitem[{Bertsimas et~al.(2018)Bertsimas, Sim, and
  Zhang}]{bertsimas2018adaptive}
Bertsimas D, Sim M, Zhang M (2018) Adaptive distributionally robust
  optimization. Management Science 65(2):604--618

\bibitem[{Birge and Louveaux(2011)}]{birge2011introduction}
Birge JR, Louveaux F (2011) Introduction to Stochastic Programming (2nd
  {E}dition). Springer Science \& Business Media

\bibitem[{Blanchet and Murthy(2019)}]{blanchet2019quantifying}
Blanchet J, Murthy K (2019) Quantifying distributional model risk via optimal
  transport. Mathematics of Operations Research 44(2):565--600

\bibitem[{Delage and Ye(2010)}]{delage2010distributionally}
Delage E, Ye Y (2010) Distributionally robust optimization under moment
  uncertainty with application to data-driven problems. Operations Research
  58(3):595--612

\bibitem[{Esfahani and Kuhn(2018)}]{esfahani2018data}
Esfahani PM, Kuhn D (2018) Data-driven distributionally robust optimization
  using the wasserstein metric: Performance guarantees and tractable
  reformulations. Mathematical Programming 171(1-2):115--166

\bibitem[{Gao and Kleywegt(2016)}]{gao2016distributionally}
Gao R, Kleywegt AJ (2016) Distributionally robust stochastic optimization with
  {W}asserstein distance. arXiv preprint arXiv:160402199

\bibitem[{Girardeau et~al.(2014)Girardeau, Leclere, and
  Philpott}]{girardeau2014convergence}
Girardeau P, Leclere V, Philpott AB (2014) On the convergence of decomposition
  methods for multistage stochastic convex programs. Mathematics of Operations
  Research 40(1):130--145

\bibitem[{Goel and Grossmann(2006)}]{goel2006class}
Goel V, Grossmann IE (2006) A class of stochastic programs with decision
  dependent uncertainty. Mathematical Programming 108(2-3):355--394

\bibitem[{Goh and Sim(2010)}]{goh2010distributionally}
Goh J, Sim M (2010) Distributionally robust optimization and its tractable
  approximations. Operations Research 58(4-part-1):902--917

\bibitem[{Guigues(2016)}]{guigues2016convergence}
Guigues V (2016) Convergence analysis of sampling-based decomposition methods
  for risk-averse multistage stochastic convex programs. SIAM Journal on
  Optimization 26(4):2468--2494

\bibitem[{Hu et~al.(2019)Hu, Li, and Mehrotra}]{hu2019data}
Hu J, Li J, Mehrotra S (2019) A data-driven functionally robust approach for
  simultaneous pricing and order quantity decisions with unknown demand
  function. Operations Research 67(6):1564--1585

\bibitem[{Jiang and Guan(2016)}]{jiang2016data}
Jiang R, Guan Y (2016) Data-driven chance constrained stochastic program.
  Mathematical Programming 158(1-2):291--327

\bibitem[{Jiang and Guan(2018)}]{jiang2018risk}
Jiang R, Guan Y (2018) Risk-averse two-stage stochastic program with
  distributional ambiguity. Operations Research 66(5):1390--1405

\bibitem[{Jonsbr{\aa}ten et~al.(1998)Jonsbr{\aa}ten, Wets, and
  Woodruff}]{jonsbraaten1998class}
Jonsbr{\aa}ten TW, Wets RJ, Woodruff DL (1998) A class of stochastic programs
  with decision dependent random elements. Annals of Operations Research
  82:83--106

\bibitem[{Kung and Liao(2018)}]{kung2018approximation}
Kung LC, Liao WH (2018) An approximation algorithm for a competitive facility
  location problem with network effects. European Journal of Operational
  Research 267(1):176--186

\bibitem[{Lappas and Gounaris(2017)}]{lappas2017use}
Lappas NH, Gounaris CE (2017) The use of decision-dependent uncertainty sets in
  robust optimization. Proceedings of Foundations of Computer-Aided Process
  Operations/Chemical Process Control 2017

\bibitem[{Lappas and Gounaris(2018)}]{lappas2018robust}
Lappas NH, Gounaris CE (2018) Robust optimization for decision-making under
  endogenous uncertainty. Computers \& Chemical Engineering 111:252--266

\bibitem[{Lee et~al.(2012)Lee, Homem-de Mello, and
  Kleywegt}]{lee2012newsvendor}
Lee S, Homem-de Mello T, Kleywegt AJ (2012) Newsvendor-type models with
  decision-dependent uncertainty. Mathematical Methods of Operations Research
  76(2):189--221

\bibitem[{Lofberg(2004)}]{lofberg2004yalmip}
Lofberg J (2004) {YALMIP}: {A} toolbox for modeling and optimization in
  {MATLAB}. In: 2004 IEEE International Conference on Robotics and Automation
  (IEEE Cat. No. 04CH37508), IEEE, pp 284--289

\bibitem[{Luo and Mehrotra(2020)}]{luo2020distributionally}
Luo F, Mehrotra S (2020) Distributionally robust optimization with decision
  dependent ambiguity sets. Optimization Letters pp 1--30

\bibitem[{McCormick(1976)}]{mccormick1976computability}
McCormick GP (1976) Computability of global solutions to factorable nonconvex
  programs: {Part I--Convex underestimating problems}. Mathematical Programming
  10(1):147--175

\bibitem[{Mehrotra and Papp(2014)}]{mehrotra2014cutting}
Mehrotra S, Papp D (2014) A cutting surface algorithm for semi-infinite convex
  programming with an application to moment robust optimization. SIAM Journal
  on Optimization 24(4):1670--1697

\bibitem[{Nohadani and Sharma(2018)}]{nohadani2018optimization}
Nohadani O, Sharma K (2018) Optimization under decision-dependent uncertainty.
  SIAM Journal on Optimization 28(2):1773--1795

\bibitem[{Noyan et~al.(2018)Noyan, Rudolf, and
  Lejeune}]{noyan2018distributionally}
Noyan N, Rudolf G, Lejeune M (2018) Distributionally robust optimization with
  decision-dependent ambiguity set. \url{http://www. optimization-online.
  org/DB HTML/2018/09/6821. html}

\bibitem[{Pereira and Pinto(1991)}]{pereira1991multi}
Pereira MV, Pinto LM (1991) Multi-stage stochastic optimization applied to
  energy planning. Mathematical Programming 52(1-3):359--375

\bibitem[{Philpott et~al.(2018)Philpott, de~Matos, and
  Kapelevich}]{philpott2018distributionally}
Philpott A, de~Matos V, Kapelevich L (2018) Distributionally robust {SDDP}.
  Computational Management Science 15(3-4):431--454

\bibitem[{Philpott and Guan(2008)}]{philpott2008convergence}
Philpott AB, Guan Z (2008) On the convergence of stochastic dual dynamic
  programming and related methods. Operations Research Letters 36(4):450--455

\bibitem[{Poss(2013)}]{poss2013robust}
Poss M (2013) Robust combinatorial optimization with variable budgeted
  uncertainty. 4OR 11(1):75--92

\bibitem[{Rockafellar and Uryasev(2002)}]{rockafellar2002conditional}
Rockafellar RT, Uryasev S (2002) Conditional value-at-risk for general loss
  distributions. Journal of Banking \& Finance 26(7):1443--1471

\bibitem[{Rockafellar et~al.(2000)Rockafellar, Uryasev
  et~al.}]{rockafellar2000optimization}
Rockafellar RT, Uryasev S, et~al. (2000) Optimization of conditional
  value-at-risk. Journal of Risk 2(3):21--42

\bibitem[{Shapiro(2001)}]{shapiro2001duality}
Shapiro A (2001) On duality theory of conic linear problems. In: Semi-infinite
  programming, Springer, pp 135--165

\bibitem[{Shapiro et~al.(2009)Shapiro, Dentcheva, and
  Ruszczy{\'n}ski}]{shapiro2009lectures}
Shapiro A, Dentcheva D, Ruszczy{\'n}ski A (2009) Lectures on Stochastic
  Programming: {M}odeling and {T}heory. SIAM

\bibitem[{Spacey et~al.(2012)Spacey, Wiesemann, Kuhn, and
  Luk}]{spacey2012robust}
Spacey SA, Wiesemann W, Kuhn D, Luk W (2012) Robust software partitioning with
  multiple instantiation. INFORMS Journal on Computing 24(3):500--515

\bibitem[{Vayanos et~al.(2011)Vayanos, Kuhn, and Rustem}]{vayanos2011decision}
Vayanos P, Kuhn D, Rustem B (2011) Decision rules for information discovery in
  multi-stage stochastic programming. In: 2011 50th IEEE Conference on Decision
  and Control and European Control Conference, IEEE, pp 7368--7373

\bibitem[{Vayanos et~al.(2020)Vayanos, Georghiou, and Yu}]{vayanos2020robust}
Vayanos P, Georghiou A, Yu H (2020) Robust optimization with decision-dependent
  information discovery. arXiv preprint arXiv:200408490

\bibitem[{Wagner(2008)}]{wagner2008stochastic}
Wagner MR (2008) Stochastic 0--1 linear programming under limited
  distributional information. Operations Research Letters 36(2):150--156

\bibitem[{Webster et~al.(2012)Webster, Santen, and
  Parpas}]{webster2012approximate}
Webster M, Santen N, Parpas P (2012) An approximate dynamic programming
  framework for modeling global climate policy under decision-dependent
  uncertainty. Computational Management Science 9(3):339--362

\bibitem[{Yu et~al.(2019)Yu, Ahmed, and Shen}]{xian2019value}
Yu X, Ahmed S, Shen S (2019) On the value of multistage stochastic facility
  location with (or without) risk aversion. Tech. rep., University of Michigan,
  Department of Industrial and Operations Engineering

\bibitem[{Zhang et~al.(2018)Zhang, Jiang, and Shen}]{zhang2018ambiguous}
Zhang Y, Jiang R, Shen S (2018) Ambiguous chance-constrained binary programs
  under mean-covariance information. SIAM Journal on Optimization
  28(4):2922--2944

\bibitem[{Zou et~al.(2019)Zou, Ahmed, and Sun}]{zou2019stochastic}
Zou J, Ahmed S, Sun XA (2019) Stochastic dual dynamic integer programming.
  Mathematical Programming 175(1-2):461--502

\end{thebibliography}


\appendix
\renewcommand{\thetheorem}{\thesection.\arabic{theorem}}
\setcounter{section}{0}
\setcounter{theorem}{0}
\renewcommand{\theequation}{\thesection-\arabic{equation}}
\setcounter{equation}{0}
\setcounter{figure}{0}
\setcounter{table}{0}

\section*{APPENDIX}

{
\section{Reformulations of N-DDDR having Continuous Supports}
\label{appen:continuous}

\normalsize
A continuous version of the Type 1 ambiguity set $\mathcal{P}_{t+1}^{D_1}(\boldsymbol{x}_{t})$ in Section \ref{sec:ambiguity_1} is given by
\begin{align}
\mathcal{P}_{t+1}^{C_1} (\boldsymbol{x}_{t}):=\left\{P\in \mathcal{M}(\Xi_{t+1},\mathcal{F}_{t+1})\ |\ \underline{\nu}(\boldsymbol{x}_t) \le P\le \bar{\nu}(\boldsymbol{x}_t), \int_{{\Xi}_{t+1}}f_s(\boldsymbol{\xi}_{t+1})P(d\boldsymbol{\xi}_{t+1})\in[l_s(\boldsymbol{x}_{t}),u_s(\boldsymbol{x}_{t})],\ \forall s\in[m]\right\},\ \label{eq:ambi_C}
\end{align}
where $\mathcal{M}(\Xi_{t+1},\mathcal{F}_{t+1})$ represents the set of all positive measures defined on $(\Xi_{t+1},\mathcal{F}_{t+1})$, and $\underline{\nu}(\boldsymbol{x}_t),\ \bar{\nu}(\boldsymbol{x}_t)\in\mathcal{M}(\Xi_{t+1},\mathcal{F}_{t+1})$ are two given measures that are lower and upper bounds for the true probability measure, respectively.  To ensure that $P$ is a probability distribution, let $l_1(\boldsymbol{x}_t)=u_1(\boldsymbol{x}_t)=f_1(\boldsymbol{\xi}_{t+1})=1$ (see \eqref{eq:normal} for more details).  

Let $\Xi_{t+1}$ be a closed and bounded set in the Euclidean space, and the probability measures $P,\ \underline{\nu},\ \bar{\nu}$ be defined on the measurable space $(\Xi_{t+1},\mathcal{F}_{t+1})$, where the $\sigma$-algebra $\mathcal{F}_{t+1}$ contains all singleton subsets, i.e., $\{\xi\}\in\mathcal{F}_{t+1}$ for all $\xi\in\Xi_{t+1}$. 
Then, based on the ambiguity set $ \mathcal{P}_{t+1}^{C_1} (\boldsymbol{x}_{t})$ in \eqref{eq:ambi_C}, we describe a reformulation of the Bellman equation \eqref{eq:bell_neutral} as an analogy to Theorem \ref{theorem:1} for the continuous support case. 
\begin{theorem}\label{theorem:1-C}
If for any feasible $\boldsymbol{x}_t\in\hat{X}_t$, the ambiguity set \eqref{eq:ambi_C} has a non-empty relative interior, then the Bellman equation \eqref{eq:bell_neutral} can be reformulated as $Q_t(\boldsymbol{x}_{t-1},\boldsymbol{\xi}_t)=$
 \small
 \begin{subequations}\label{Bellman_C}
\begin{align*}
\min_{\boldsymbol{\alpha},\boldsymbol{\beta},{\underline{\gamma},\bar{\gamma},}\boldsymbol{x}_t,\boldsymbol{y}_t}\quad &g_t(\boldsymbol{x}_t,\boldsymbol{y}_t)-{\boldsymbol{\alpha}}^{\mathsf T} \boldsymbol{l}(\boldsymbol{x}_t)+{\boldsymbol{\beta}}^{\mathsf T} \boldsymbol{u}(\boldsymbol{x}_t) {-\int_{\boldsymbol{\xi}_{t+1}\in\Xi_{t+1}}\underline{\gamma}(\boldsymbol{\xi}_{t+1})\underline{\nu}(\boldsymbol{x}_t,\boldsymbol{\xi}_{t+1})d\boldsymbol{\xi}_{t+1} + \int_{\boldsymbol{\xi}_{t+1}\in\Xi_{t+1}}\bar{\gamma}(\boldsymbol{\xi}_{t+1})\bar{\nu}(\boldsymbol{x}_t,\boldsymbol{\xi}_{t+1})d\boldsymbol{\xi}_{t+1}}\\
\text{s.t.}\quad&(-\boldsymbol{\alpha}+\boldsymbol{\beta})^{\mathsf T} \boldsymbol{f}(\boldsymbol{\xi}_{t+1}) {-\underline{\gamma}(\boldsymbol{\xi}_{t+1})+\bar{\gamma}(\boldsymbol{\xi}_{t+1})}\ge Q_{t+1}(\boldsymbol{x}_t,\boldsymbol{\xi}_{t+1}),\ \forall \boldsymbol{\xi}_{t+1}\in \Xi_{t+1},\\
&(\boldsymbol{x}_t,\boldsymbol{y}_t)\in X_t(\boldsymbol{x}_{t-1},\boldsymbol{\xi}_t),\\
& \boldsymbol{\alpha},\ \boldsymbol{\beta},\  {\boldsymbol{\underline{\gamma}}(\boldsymbol{\xi}_{t+1}),\ \boldsymbol{\bar{\gamma}}(\boldsymbol{\xi}_{t+1})}\ge 0, \ \forall \boldsymbol{\xi}_{t+1}\in\Xi_{t+1}.
\end{align*}
 \end{subequations}
 \end{theorem}
   \normalsize 
 
Next, a continuous version of the Type 2 ambiguity set $\mathcal{P}_{t+1}^{D_2} (\boldsymbol{x}_{t})$ in \eqref{eq:ambi_match} is given by
\begin{subequations}
\label{eq:ambi_match_C}
  \begin{align}
\mathcal{P}^{C_2}_{t+1}(\boldsymbol{x}_t):=\Biggl \{P\in \mathcal{M}(\Xi_{t+1},\mathcal{F}_{t+1})\ |\  &\mathbb{E}_P[\boldsymbol{\xi}_{t+1}]=\boldsymbol{\mu}(\boldsymbol{x}_t),\label{eq:ambi_match1}\\
& \mathbb{E}_P[(\boldsymbol{\xi}_{t+1}-\boldsymbol{\mu}(\boldsymbol{x}_t))(\boldsymbol{\xi}_{t+1}-\boldsymbol{\mu}(\boldsymbol{x}_t))^{\mathsf T}]=\boldsymbol{\Sigma}(\boldsymbol{x}_t)\Biggr \},\label{eq:ambi_match2}
\end{align}
\end{subequations}
Then, the following result is an analogy to Theorem \ref{theorem:2} based on the continuous ambiguity set \eqref{eq:ambi_match_C}. 
\begin{theorem}
\label{theorem:2-C}
If for any feasible $\boldsymbol{x}_t\in\hat{X}_t$, the ambiguity set \eqref{eq:ambi_match_C} has a non-empty relative interior, then the Bellman equation \eqref{eq:bell_neutral} can be reformulated as
\small
\begin{subequations}\label{Bellman_2_C}
 \begin{align}
Q_t(\boldsymbol{x}_{t-1},\boldsymbol{\xi}_t)=\min_{\boldsymbol{x}_t,\boldsymbol{y}_t,s,\boldsymbol{u},\boldsymbol{Y}}\quad &g_t(\boldsymbol{x}_t,\boldsymbol{y}_t)+s+\boldsymbol{u}^{\mathsf T}\boldsymbol{\mu}(\boldsymbol{x}_t)+\boldsymbol{\Sigma}(\boldsymbol{x}_t)\bullet \boldsymbol{Y}\\
\text{s.t.}\quad&s+\boldsymbol{u}^{\mathsf T}\boldsymbol{\xi}_{t+1} +(\boldsymbol{\xi}_{t+1}-\boldsymbol{\mu}(\boldsymbol{x}_t))(\boldsymbol{\xi}_{t+1}-\boldsymbol{\mu}(\boldsymbol{x}_t))^{\mathsf T}\bullet \boldsymbol{Y}\ge Q_{t+1}(\boldsymbol{x}_t,\boldsymbol{\xi_{t+1}}),\nonumber\\
&\hspace{12ex}\forall \boldsymbol{\xi}_{t+1}\in\Xi_{t+1},\\
&(\boldsymbol{x}_t,\boldsymbol{y}_t)\in X_t(\boldsymbol{x}_{t-1},\boldsymbol{\xi}_t).\nonumber
\end{align}
\end{subequations}
\end{theorem}
\normalsize

Finally, a continuous version of the Type 3 ambiguity set $ \mathcal{P}_{t+1}^{D_3}(\boldsymbol{x}_{t})$ in \eqref{eq:ambi_bound} is given by
\begin{subequations}
\label{eq:ambi_bound_C}
 \begin{align}
\mathcal{P}_{t+1}^{C_3}(\boldsymbol{x}_t):=\Biggl \{P\in \mathcal{M}(\Xi_{t+1},\mathcal{F}_{t+1})\ |\ 
&(\mathbb{E}_P[\boldsymbol{\xi}_{t+1}]-\boldsymbol{\mu}(\boldsymbol{x}_t))^{\mathsf T}\boldsymbol{\Sigma}(\boldsymbol{x}_t)^{-1}(\mathbb{E}_P[\boldsymbol{\xi}_{t+1}]-\boldsymbol{\mu}(\boldsymbol{x}_t))\le \gamma,\\
& \mathbb{E}_P[(\boldsymbol{\xi}_{t+1}-\boldsymbol{\mu}(\boldsymbol{x}_t))(\boldsymbol{\xi}_{t+1}-\boldsymbol{\mu}(\boldsymbol{x}_t))^{\mathsf T}]\preceq \eta \boldsymbol{\Sigma}(\boldsymbol{x}_t)\Biggr \}. 
\end{align}
\end{subequations}
\normalsize
We present a reformulation of  the Bellman equation \eqref{eq:bell_neutral} in the following theorem that is an analogy to Theorem \ref{theorem:3}, but given the continuous ambiguity set \eqref{eq:ambi_bound_C}.
\begin{theorem}
\label{theorem:3-C}
Suppose that the Slater's constraint qualification conditions
are satisfied, i.e., for any feasible $\boldsymbol{x}_t\in \hat{X}_t$, there exists a probability measure $P\in \mathcal{M}(\Xi_{t+1},\mathcal{F}_{t+1})$ such that $(\mathbb{E}_P[\boldsymbol{\xi}_{t+1}]-\boldsymbol{\mu}(\boldsymbol{x}_t))^{\mathsf T}\boldsymbol{\Sigma}(\boldsymbol{x}_t)^{-1}(\mathbb{E}_P[\boldsymbol{\xi}_{t+1}]-\boldsymbol{\mu}(\boldsymbol{x}_t))< \gamma$, and $\mathbb{E}_P[(\boldsymbol{\xi}_{t+1}-\boldsymbol{\mu}(\boldsymbol{x}_t))(\boldsymbol{\xi}_{t+1}-\boldsymbol{\mu}(\boldsymbol{x}_t))^{\mathsf T}]\prec \eta \boldsymbol{\Sigma}(\boldsymbol{x}_t)$.
Using the ambiguity set defined in \eqref{eq:ambi_bound_C}, the Bellman equation \eqref{eq:bell_neutral} can be recast as
\small
\begin{subequations}\label{eq:bell_bound1_C}
 \begin{align}
Q_t(\boldsymbol{x}_{t-1},\boldsymbol{\xi}_t)=\min_{\boldsymbol{x}_t,\boldsymbol{y}_t,s,\boldsymbol{Z},\boldsymbol{Y}}\quad &g_t(\boldsymbol{x}_t,\boldsymbol{y}_t)+s+\boldsymbol{\Sigma}(\boldsymbol{x}_t)\bullet \boldsymbol{z}_1 -2\boldsymbol{\mu}(\boldsymbol{x}_t)^{\mathsf T}\boldsymbol{z}_2+\gamma z_3 +\eta\boldsymbol{\Sigma}(\boldsymbol{x}_t)\bullet \boldsymbol{Y}\\
\text{s.t.}\quad&s-2\boldsymbol{z}_2^{\mathsf T}\boldsymbol{\xi}_{t+1}+(\boldsymbol{\xi}_{t+1}-\boldsymbol{\mu}(\boldsymbol{x}_t))(\boldsymbol{\xi}_{t+1}-\boldsymbol{\mu}(\boldsymbol{x}_t))^{\mathsf T}\bullet Y\ge Q_{t+1}(\boldsymbol{x}_t,\boldsymbol{\xi}_{t+1}),\nonumber\\
& \hspace{12ex}\forall \boldsymbol{\xi}_{t+1}\in \Xi_{t+1},\\
& \boldsymbol{Z}=\begin{pmatrix} \boldsymbol{z}_1&\boldsymbol{z}_2\\\boldsymbol{z}_2^{\mathsf T}&z_3\end{pmatrix}\succeq 0,\ \boldsymbol{Y}\succeq 0,\\
&(\boldsymbol{x}_t,\boldsymbol{y}_t)\in X_t(\boldsymbol{x}_{t-1},\boldsymbol{\xi}_t).
\end{align}
\end{subequations}
\end{theorem}
\normalsize
The proof of Theorem \ref{theorem:1-C} is similar to the proof in \citep{luo2020distributionally} for the two-stage continuous-support case, and we omit its details. We provide detailed proofs for Theorems \ref{theorem:2-C} and \ref{theorem:3-C} in Appendix \ref{appen:allproofs}. Note that the above three reformulations for ambiguity sets with continuous support are semi-infinite programs and thus cannot be optimized directly.
}

\section{Risk-averse Multistage DRO with Endogenous Uncertainty}
\label{sec:dro_risk}
\normalsize
We can extend N-DDDR in \eqref{eq:1} to a more general setting. Previously, the robust counterpart chooses the worst-case distribution $P$ from a risk-neutral aspect using expectation to measure the uncertain cost over multiple stages. However, a decision maker may measure the worst-case distribution in a risk-averse fashion, and we accordingly replace the expectations by coherent risk measures $\rho_t, \ \forall t=2,\ldots,T$. The corresponding risk-averse multistage decision-dependent DRO model is: 
\begin{align}
\mbox{{\bf A-DDDR:}} & ~\nonumber\\
\quad\min_{(\boldsymbol{\boldsymbol{x}_1},\boldsymbol{y}_1)\in X_1}&\Big\{g_1(\boldsymbol{x}_1,\boldsymbol{y}_1)+\max_{P_2\in \mathcal{P}_2(\boldsymbol{x}_1)}\rho_2\Big[\min_{(\boldsymbol{x}_2,\boldsymbol{y}_2)\in X_2(\boldsymbol{x}_1,\boldsymbol{\xi}_2)}g_2(\boldsymbol{x}_2,\boldsymbol{y}_2)+\cdots\nonumber\\
&+\max_{P_{t}\in\mathcal{P}_{t}(\boldsymbol{x}_{t-1})}\rho_{t}\Big[\min_{(\boldsymbol{x}_{t},\boldsymbol{y}_{t})\in X_{t}(\boldsymbol{x}_{t-1},\boldsymbol{\xi}_{t})}g_{t}(\boldsymbol{x}_{t},\boldsymbol{y}_{t})+\cdots\nonumber\\
&+\max_{P_T\in\mathcal{P}_T(\boldsymbol{x}_{T-1})}\rho_{T}\Big[\min_{(\boldsymbol{x}_T,\boldsymbol{y}_T)\in X_T(\boldsymbol{x}_{T-1},\boldsymbol{\xi}_T)}g_T(\boldsymbol{x}_T,\boldsymbol{y}_T)\Big]\Big\}.\label{eq:risk}
\end{align}
We consider a special class of coherent risk measures, which is a convex combination of expectation and Conditional Value-at-Risk (CVaR) \citep{rockafellar2000optimization}:
\begin{equation*}
\rho_t(Z)=(1-\lambda_t)\mathbb{E}[Z]+\lambda_t \text{CVaR}_{\alpha_t}[Z],
\label{eq:rho}
\end{equation*}
where $\lambda_t\in[0,1]$ is a parameter that balances the expectation and CVaR measure at $\alpha_t\in(0,1)$ risk level. This risk measure is more general than expectation and it becomes the risk-neutral case when $\lambda_t=0$. 

The Bellman equations for A-DDDR \eqref{eq:risk} then become:
\begin{equation*}
Q_1=\min_{(\boldsymbol{x}_1,\boldsymbol{y}_1)\in X_1} g_1(\boldsymbol{x}_1,\boldsymbol{y}_1)+\max_{P_2\in\mathcal{P}_2(\boldsymbol{x}_1)}\rho_2[Q_2(\boldsymbol{x}_1,\boldsymbol{\xi}_2)],
\end{equation*}
where for $t=2,\ldots,T-1$,
\begin{equation}
Q_t(\boldsymbol{x}_{t-1},\boldsymbol{\xi}_t)=\min_{(\boldsymbol{x}_t,\boldsymbol{y}_t)\in X_t(\boldsymbol{x}_{t-1},\boldsymbol{\xi}_t)} g_t(\boldsymbol{x}_t,\boldsymbol{y}_t)+\max_{P_{t+1}\in\mathcal{P}_{t+1}(\boldsymbol{x}_t)}\rho_{t+1}[Q_{t+1}(\boldsymbol{x}_t,\boldsymbol{\xi}_{t+1})],\label{eq:bell_risk}
\end{equation}
and
\begin{equation*}
Q_T(\boldsymbol{x}_{T-1},\boldsymbol{\xi}_T)=\min_{(\boldsymbol{x}_T,\boldsymbol{y}_T)\in X_T(\boldsymbol{x}_{T-1},\boldsymbol{\xi}_T)} g_T(\boldsymbol{x}_T,\boldsymbol{y}_T).\label{eq:bell_risk_T}
\end{equation*}

Following the results by \citet{rockafellar2002conditional}, CVaR can be attained by solving the following optimization problem:
\begin{equation*}
\text{CVaR}_{\alpha_t}[Z]:=\inf_{\eta\in\mathbb{R}}\left\lbrace \eta+\frac{1}{1-\alpha_t}\mathbb{E}[Z-\eta]_{+}\right\rbrace ,\label{eq:cvar}
\end{equation*}
where $[a]_{+}:=\max\{a,0\},$ and $\eta$ is an auxiliary variable. To linearize $[Z-\eta]_{+}$, we replace it by a variable $m$ with two additional constraints: $m\ge 0,\ m\ge Z-\eta$.

Recall that by assumption, every $P_{t+1}\in\mathcal{P}_{t+1}(\boldsymbol{x}_{t})$ has a decision-independent finite support $\Xi_{t+1}:=\{\boldsymbol{\xi}_{t+1}^k\}_{k=1}^K,\ \forall \boldsymbol{x}_{t}\in X_{t}$ for a fixed $K$ and all $t\in[T-1]$. Each realization $k\in[K]$ is associated with probability $p_k$, and therefore the inner maximization problem $\max_{P_{t+1}\in\mathcal{P}_{t+1}(\boldsymbol{x}_t)}\rho_{t+1}[Q_{t+1}(\boldsymbol{x}_t,\boldsymbol{\xi}_{t+1})]$ in \eqref{eq:bell_risk} can be reformulated as 
\begin{subequations}\label{eq:inner2}
\begin{align}
\max_{P_{t+1}\in\mathcal{P}_{t+1}(\boldsymbol{x}_t)}\min_{m,\eta}\quad & \lambda_{t+1}\eta+\sum_{k=1}^Kp_k(\frac{\lambda_{t+1}}{1-\alpha_{t+1}}m+(1-\lambda_{t+1})Q_{t+1}^k)\\
\text{s.t.}\quad&m+\eta\ge Q_{t+1}^k,\ \forall k\in[K],\label{con}\\
& m\ge 0,\ \forall k\in[K].
\end{align}
\end{subequations}
We further simplify the notation of the recursive function $Q_{t+1}(\boldsymbol{x}_t,\boldsymbol{\xi}_{t+1}^k)$ as $Q_{t+1}^k$. 
Associating dual variables $q_k$ with constraints \eqref{con} and applying strong duality result, we have
\small
\begin{subequations}
\begin{align}
    \max_{P_{t+1}\in\mathcal{P}_{t+1}(\boldsymbol{x}_t)}\Biggl \{(1-\lambda_{t+1})\sum_{k=1}^Kp_kQ_{t+1}^k+\max_{\boldsymbol{q}}\quad &\sum_{k=1}^Kq_kQ_{t+1}^k\\
    \text{s.t.}\quad& q_k\le p_k\frac{\lambda_{t+1}}{1-\alpha_{t+1}},\ \forall k\in[K],\label{eq:w1}\\
    & \sum_{k=1}^Kq_k=\lambda_{t+1},\label{eq:w2}\\
    &q_k\ge 0,\ \forall k\in[K].\Biggr \}
\end{align}
\end{subequations}
\normalsize
Merging the two layers of maximization problems, for each $t\in[T-1]$, we solve
\begin{subequations}\label{eq:riskbell}
\begin{align}
    \max_{\boldsymbol{p,q}}\quad&(1-\lambda_{t+1})\sum_{k=1}^K p_kQ_{t+1}^k+q_kQ_{t+1}^k\\
    \text{s.t.}\quad& \mbox{\eqref{eq:w1},\ \eqref{eq:w2}},\nonumber\\
    & \boldsymbol{p}\in \mathcal{P}_{t+1}(\boldsymbol{x}_t),\\
    &q_k\ge 0,\ \forall k\in[K] .
\end{align}
\end{subequations}
In the following subsections, we present reformulations of A-DDDR in \eqref{eq:risk} under the three types of ambiguity sets mentioned in Section \ref{sec:dro}.

\subsection{Solving A-DDDR under Type 1 Ambiguity Set}
Using the ambiguity set defined in \eqref{eq:ambi}, the inner maximization problem \eqref{eq:riskbell} can be recast as 
\begin{subequations}\label{eq:feasible_3}
\begin{align}
    \max_{\boldsymbol{p,q}}\quad&(1-\lambda_{t+1})\sum_{k=1}^Kp_kQ_{t+1}^k+q_kQ_{t+1}^k\\
    \text{s.t.}\quad& \mbox{\eqref{eq:w1},\ \eqref{eq:w2},\ \eqref{eq:lower}--\eqref{eq:upper2}}\nonumber\\
    &p_k,\ q_k\ge 0,\ \forall k\in[K] .
\end{align}
\end{subequations}

\begin{theorem}
\label{theorem:4}
If for any feasible $\boldsymbol{x}_t\in \hat{X}_t$, problem \eqref{eq:feasible_3} is feasible, then the Bellman equation \eqref{eq:bell_risk} can be reformulated as $Q_t(\boldsymbol{x}_{t-1},\boldsymbol{\xi}_t)=$
\small
\begin{subequations}\label{Bellman2}
\begin{align}
\min_{\boldsymbol{\alpha},\boldsymbol{\beta},\boldsymbol{x}_t,\boldsymbol{y}_t} \quad&g_t(\boldsymbol{x}_t,\boldsymbol{y}_t)+\lambda_{t+1}\theta-{\boldsymbol{\alpha}}^{\mathsf T} \boldsymbol{l}(\boldsymbol{x}_t)+{\boldsymbol{\beta}}^{\mathsf T} \boldsymbol{u}(\boldsymbol{x}_t){-\underline{\boldsymbol{\gamma}}^{\mathsf T}\underline{\boldsymbol{p}}(\boldsymbol{x}_t) + \bar{\boldsymbol{\gamma}}^{\mathsf T}\bar{\boldsymbol{p}}(\boldsymbol{x}_t)}\\
\text{s.t.}\quad& \pi_k+\theta\ge Q_{t+1}^k,\ \forall k\in[K],\\
&-\frac{\lambda_{t+1}}{1-\alpha_{t+1}}\pi_k+(-\boldsymbol{\alpha}+\boldsymbol{\beta})^{\mathsf T} \boldsymbol{f}(\boldsymbol{\xi}_{t+1}^k){-\underline{\gamma}_k+\bar{\gamma}_k}\ge (1-\lambda_{t+1}) Q_{t+1}^k,\ \forall k\in[K],\\
&(\boldsymbol{x}_t,\boldsymbol{y}_t)\in X_t(\boldsymbol{x}_{t-1},\boldsymbol{\xi}_t),\\
&\pi_k,\ \boldsymbol{\alpha},\ \boldsymbol{\beta},\ {\boldsymbol{\underline{\gamma}},\ \boldsymbol{\bar{\gamma}}}\ge 0.
\end{align}
\end{subequations}
\normalsize
\end{theorem}
The proof of Theorem \ref{theorem:4} is similar to the one of Theorem \ref{theorem:1} in Appendix \ref{appen:allproofs}, where the only difference is that we introduce the two more dual variables $\pi_k$ and $\theta$, associated with constraints \eqref{eq:w1} and \eqref{eq:w2}, respectively. 

Notice here when $\lambda_{t+1}=0$, model \eqref{Bellman2} reduces to the risk-neutral case \eqref{Bellman}. This reformulation also has similar computational complexity as model \eqref{Bellman} in Theorem \ref{theorem:1}. Therefore, with the same specific ambiguity set considered in \eqref{eq:ambi_ineq} in Section \ref{sec:dro}, we can apply McCormick envelopes to obtain a multistage stochastic MILP and deploy SDDiP to solve it. 

\subsection{Solving A-DDDR under Type 2 Ambiguity Set}
Under Type 2 ambiguity set in \eqref{eq:ambi_match}, the inner maximization problem \eqref{eq:riskbell} can be recast as 
\begin{subequations}\label{eq:2}
\begin{align}
    \max_{\boldsymbol{p,q}}\quad&(1-\lambda_{t+1})\sum_{k=1}^K p_kQ_{t+1}^k+q_kQ_{t+1}^k\\
    \text{s.t.}\quad& \mbox{\eqref{eq:w1},\ \eqref{eq:w2},\ \eqref{eq:match_prob1}--\eqref{eq:match_cov}}\nonumber\\
& p_k,\ q_k\ge0, \quad \forall k\in[K].
\end{align}
\end{subequations}

\begin{theorem}
\label{theorem:5}
If for any feasible $\boldsymbol{x}_t\in \hat{X}_t$, problem \eqref{eq:2} is feasible, then the Bellman equation \eqref{eq:bell_risk} can be reformulated as $Q_t(\boldsymbol{x}_{t-1},\boldsymbol{\xi}_t)=$
\small
\begin{subequations}\label{eq:refor}
 \begin{align}
\min_{\boldsymbol{x}_t,\boldsymbol{y}_t,s,\boldsymbol{u},\boldsymbol{Y}} \quad&g_t(\boldsymbol{x}_t,\boldsymbol{y}_t)-\lambda_{t+1}\theta+s+\boldsymbol{u}^{\mathsf T}\boldsymbol{\mu}(\boldsymbol{x}_t)+\boldsymbol{\Sigma}(\boldsymbol{x}_t)\bullet \boldsymbol{Y} \\
\text{s.t.}\quad&s+\boldsymbol{u}^{\mathsf T}\boldsymbol{\xi}_{t+1}^k +(\boldsymbol{\xi}_{t+1}^k-\boldsymbol{\mu}(\boldsymbol{x}_t))(\boldsymbol{\xi}_{t+1}^k-\boldsymbol{\mu}(\boldsymbol{x}_t))^{\mathsf T}\bullet \boldsymbol{Y}-\pi_k\frac{\lambda_{t+1}}{1-\alpha_{t+1}}\ge (1-\lambda_{t+1})Q_{t+1}^k,\nonumber\\
& \hspace{15ex} \forall k\in[K],\\
& \pi_k-\theta\ge Q_{t+1}^k,\quad \forall k\in[K],\\
& \pi_k\ge 0 ,\quad \forall k\in[K],\\
&(\boldsymbol{x}_t,\boldsymbol{y}_t)\in X_t(\boldsymbol{x}_{t-1},\boldsymbol{\xi}_t).
\end{align}
\end{subequations}
\normalsize
\end{theorem}
The proof of Theorem \ref{theorem:5} is similar to the one of Theorem \ref{theorem:2} in Appendix \ref{appen:allproofs}, where the only difference is that we introduce the two more dual variables $\pi_k$ and $\theta$, associated with constraints \eqref{eq:w1} and \eqref{eq:w2}, respectively. 

Notice here when $\lambda_{t+1}=0$, the risk-averse model \eqref{eq:refor} reduces to the risk-neutral case \eqref{Bellman_2}. We can apply McCormick envelopes to get a multistage stochastic MILP and use SDDiP algorithm to attain optimal solutions as in Section \ref{sec:dro_match}.

\subsection{	Solving A-DDDR under Type 3 Ambiguity Set}
Given Type 3 ambiguity set defined in \eqref{eq:ambi_bound}, the inner maximization problem \eqref{eq:riskbell} can be recast as 
\begin{align*}
    \max_{\boldsymbol{p,q}}\quad&(1-\lambda_{t+1})\sum_{k=1}^Kp_kQ_{t+1}^k+q_kQ_{t+1}^k\\
    \text{s.t.}\quad& \mbox{\eqref{eq:w1},\ \eqref{eq:w2},\ \eqref{eq:equal}--\eqref{eq:Sigma}}\nonumber\\
    & p_k\ge 0,\ q_k\ge0, \quad \forall k\in[K].
\end{align*}

\begin{theorem}
\label{theorem:6}
Suppose that Slater's constraint qualification conditions
are satisfied, i.e., for any feasible $\boldsymbol{x}_t\in \hat{X}_t$, there exists a vector $p = (p_1,p_2,\ldots,p_K)^{\mathsf T}$ such that $\sum_{k=1}^Kp_{k}=1$,  $(\sum_{k=1}^K p_k\boldsymbol{\xi}_{t+1}^k-\boldsymbol{\mu}(\boldsymbol{x}_t))^{\mathsf T}\boldsymbol{\Sigma}(\boldsymbol{x}_t)^{-1}(\sum_{k=1}^K p_k\boldsymbol{\xi}_{t+1}^k-\boldsymbol{\mu}(\boldsymbol{x}_t))< \gamma$, and $\sum_{k=1}^Kp_k(\boldsymbol{\xi}_{t+1}^k-\boldsymbol{\mu}(\boldsymbol{x}_t))(\boldsymbol{\xi}_{t+1}^k-\boldsymbol{\mu}(\boldsymbol{x}_t))^{\mathsf T}\prec \eta \boldsymbol{\Sigma}(\boldsymbol{x}_t)$.
Using the ambiguity set defined in \eqref{eq:ambi_bound}, the Bellman equation \eqref{eq:bell_risk} can be recast as $Q_t(\boldsymbol{x}_{t-1},\boldsymbol{\xi}_t)=$
\small
\begin{subequations}\label{eq:bell_bound2}
 \begin{align}
\min_{\boldsymbol{x}_t,\boldsymbol{y}_t,s,\boldsymbol{Z},\boldsymbol{Y}}\quad &g_t(\boldsymbol{x}_t,\boldsymbol{y}_t)-\lambda\theta+s+\boldsymbol{\Sigma}(\boldsymbol{x}_t)\bullet \boldsymbol{z}_1 -2\boldsymbol{\mu}(\boldsymbol{x}_t)^{\mathsf T}\boldsymbol{z}_2+\gamma z_3 +\eta\boldsymbol{\Sigma}(\boldsymbol{x}_t)\bullet \boldsymbol{Y}\\
\text{s.t.}\quad&s-2\boldsymbol{z}_2^{\mathsf T}\boldsymbol{\xi}_{t+1}^k+(\boldsymbol{\xi}_{t+1}^k-\boldsymbol{\mu}(\boldsymbol{x}_t))(\boldsymbol{\xi}_{t+1}^k-\boldsymbol{\mu}(\boldsymbol{x}_t))^{\mathsf T}\bullet Y-\pi_k\frac{\lambda_{t+1}}{1-\alpha_{t+1}}\ge (1-\lambda_{t+1})Q_{t+1}^k,\nonumber\\
&\hspace{48ex}\forall k\in[K],\\
& \boldsymbol{Z}=\begin{pmatrix} \boldsymbol{z}_1&\boldsymbol{z}_2\\\boldsymbol{z}_2^{\mathsf T}&z_3\end{pmatrix}\succeq 0,\ \boldsymbol{Y}\succeq 0,\\
& \pi_k - \theta \ge Q_{t+1}^k,\ \pi_k\ge 0 ,\quad \forall k\in[K],\\
&(\boldsymbol{x}_t,\boldsymbol{y}_t)\in X_t(\boldsymbol{x}_{t-1},\boldsymbol{\xi}_t).
\end{align}
\end{subequations}
\normalsize
\end{theorem}
The proof is similar to the one of Theorem \ref{theorem:3}. All the proofs in this section are omitted here due to similarity. 

Notice at we obtain an MISDP in each stage and when $\lambda_{t+1}=0$, the risk-averse model \eqref{eq:bell_bound2} reduces to the risk-neutral case \eqref{eq:bell_bound1}.  We can apply McCormick envelopes and approximation schemes to obtain valid upper and lower bounds similar to the procedures in Sections \ref{sec:relax}--\ref{sec:upper}.

\section{Details of All Needed Proofs}
\label{appen:allproofs}

\begin{proof}[Theorem \ref{theorem:1}]
The proof follows Theorem 3.1 in \citep{luo2020distributionally}.
Using the ambiguity set defined in \eqref{eq:ambi}, the inner maximization problem of \eqref{eq:bell_neutral} can be expressed as
\begin{subequations}\label{eq:inner}
\begin{align}
\max_{\boldsymbol{p}\in\mathbb{R}^K}\quad&\sum_{k=1}^Kp_kQ_{t+1}(\boldsymbol{x}_t,\boldsymbol{\xi}_{t+1}^k)\\
\text{s.t.}\quad&\sum_{k=1}^Kp_k\boldsymbol{f}(\boldsymbol{\xi}_{t+1}^k)\ge\boldsymbol{l}(\boldsymbol{x}_{t}),\label{eq:lower}\\
&\sum_{k=1}^Kp_k\boldsymbol{f}(\boldsymbol{\xi}_{t+1}^k)\le\boldsymbol{u}(\boldsymbol{x}_t),\label{eq:upper}\\
&{p_k\ge \underline{p}_k(\boldsymbol{x}_t),\ \forall k\in [K],}\label{eq:lower2}\\
& {p_k\le \bar{p}_k(\boldsymbol{x}_t),\ \forall k\in [K],} \label{eq:upper2}\\
& p_k\ge 0, \ \forall k\in[K].
\end{align}
\end{subequations}
We associate dual variables $\boldsymbol{\alpha},\ \boldsymbol{\beta}\in \mathbb{R}^m$ with Constraints \eqref{eq:lower} and \eqref{eq:upper}, {dual variables $\boldsymbol{\underline{\gamma}}$ and $\boldsymbol{\bar{\gamma}}\in \mathbb{R}^K$ with Constraints \eqref{eq:lower2} and \eqref{eq:upper2}}, respectively. When \eqref{eq:inner} is feasible, 
strong duality holds and the Bellman equation \eqref{eq:bell_neutral} can be reformulated as \eqref{Bellman},
which completes the proof. $\square$
\end{proof}

The proof of Theorem \ref{theorem:1-C} is omitted due to its similarity to the proof of Theorem 3.3 in \citep{luo2020distributionally}. 

\begin{proof}[Theorem \ref{theorem:2}]
Following Type 2 ambiguity set in \eqref{eq:ambi_match}, the inner maximization problem in \eqref{eq:bell_neutral} can be recast as
\begin{subequations}\label{eq:feasible_2}
\begin{align}
\max_{\boldsymbol{p}\in\mathbb{R}^K}\quad&\sum_{k=1}^Kp_kQ_{t+1}(\boldsymbol{x}_t,\boldsymbol{\xi}_{t+1}^k)\\
\text{s.t.}\quad&\sum_{k=1}^K p_k=1,\label{eq:match_prob1}\\
&\sum_{k=1}^K p_k\boldsymbol{\xi}_{t+1}^k=\boldsymbol{\mu}(\boldsymbol{x}_t), \\
&\sum_{k=1}^K p_k(\boldsymbol{\xi}_{t+1}^k-\boldsymbol{\mu}(\boldsymbol{x}_t))(\boldsymbol{\xi}_{t+1}^k-\boldsymbol{\mu}(\boldsymbol{x}_t))^{\mathsf T}=\boldsymbol{\Sigma}(\boldsymbol{x}_t), \label{eq:match_cov}\\
& p_k\ge 0,\ \forall k\in[K].
\end{align}
\end{subequations}
If the above linear program is feasible, then strong duality holds. Associate dual variables $s\in\mathbb{R},\ \boldsymbol{u}\in\mathbb{R}^J, \ \boldsymbol{Y}\in\mathbb{R}^{J\times J}$ with the three sets of constraints, respectively, and recast the inner maximization problem as a minimization problem. After including constraints $(\boldsymbol{x}_t,\boldsymbol{y}_t)\in X_t(\boldsymbol{x}_{t-1},\boldsymbol{\xi}_t)$, the Bellman equation \eqref{eq:bell_neutral} is equivalent to  \eqref{Bellman_2},
and we complete the proof. $\square$
\end{proof}

\begin{proof}[Theorem \ref{theorem:2-C}]
The proof follows the conic duality in functional spaces \citep{shapiro2001duality}. 
Using the ambiguity set defined in \eqref{eq:ambi_match_C}, the inner maximization problem of \eqref{eq:bell_neutral} can be formulated as a conic linear program in a functional space as follows:
\begin{subequations}\label{eq:inner_2_C}
\begin{align}
\max_{P\in \mathcal{M}(\Xi_{t+1},\mathcal{F}_{t+1})}\quad&\mathbb{E}_P[Q_{t+1}(\boldsymbol{x}_t,\boldsymbol{\xi}_{t+1})]\\
\text{s.t.}\quad& \mathbb{E}_P[1]= 1,\label{eq:inner_2_C_1}\\
& \mathbb{E}_P[\boldsymbol{\xi}_{t+1}]=\boldsymbol{\mu}(\boldsymbol{x}_t),\label{eq:inner_2_C_2}\\
& \mathbb{E}_P[(\boldsymbol{\xi}_{t+1}-\boldsymbol{\mu}(\boldsymbol{x}_t))(\boldsymbol{\xi}_{t+1}-\boldsymbol{\mu}(\boldsymbol{x}_t))^{\mathsf T}]=\boldsymbol{\Sigma}(\boldsymbol{x}_t) \label{eq:inner_2_C_3}
\end{align}
\end{subequations}
 We associate dual variables $s\in\mathbb{R},\ \boldsymbol{u}\in\mathbb{R}^J, \ \boldsymbol{Y}\in\mathbb{R}^{J\times J}$ with the three sets of constraints, respectively.
Because the primal problem has a non-empty relative interior, strong duality holds and the dual problem can be formulated as \eqref{Bellman_2_C}, which completes the proof. $\square$
\end{proof}

\begin{proof}[Theorem \ref{theorem:3}]
Given Type 3 ambiguity set \eqref{eq:ambi_bound}, the inner maximization problem in \eqref{eq:bell_neutral} can be recast as
\begin{subequations}
\label{eq:inner_bound}
\begin{align}
\max_{\boldsymbol{p}\in\mathbb{R}^K, \boldsymbol{\tau}\in\mathbb{R}^J}\quad&\sum_{k=1}^Kp_kQ_{t+1}(\boldsymbol{x}_t,\boldsymbol{\xi}_{t+1}^k)\\
\text{s.t.}\quad&\sum_{k=1}^K p_k=1, \label{eq:equal}\\
&\sum_{k=1}^K p_k\boldsymbol{\xi}_{t+1}^k=\boldsymbol{\tau}, \\
& (\boldsymbol{\tau}-\boldsymbol{\mu}(\boldsymbol{x}_t))^{\mathsf T}\boldsymbol{\Sigma}(\boldsymbol{x}_t)^{-1}(\boldsymbol{\tau}-\boldsymbol{\mu}(\boldsymbol{x}_t))\le \gamma, \label{eq:inner_sdp}\\
&\sum_{k=1}^K p_k(\boldsymbol{\xi}_{t+1}^k-\boldsymbol{\mu}(\boldsymbol{x}_t))(\boldsymbol{\xi}_{t+1}^k-\boldsymbol{\mu}(\boldsymbol{x}_t))^{\mathsf T}\preceq \eta \boldsymbol{\Sigma}(\boldsymbol{x}_t),\label{eq:Sigma} \\
& p_k\ge 0,\ \forall k\in[K].
\end{align}
\end{subequations}
We rewrite Constraint \eqref{eq:inner_sdp} as
\begin{equation*}
    \begin{pmatrix}
    \boldsymbol{\Sigma}(\boldsymbol{x}_t) & \boldsymbol{\tau}-\boldsymbol{\mu}(\boldsymbol{x}_t)\\
    (\boldsymbol{\tau}-\boldsymbol{\mu}(\boldsymbol{x}_t))^{\mathsf T} & \gamma
    \end{pmatrix} \succeq 0,
\end{equation*}
and associate dual variables $s\in\mathbb{R},\ \boldsymbol{u}\in\mathbb{R}^d, \ \boldsymbol{Z}=\begin{pmatrix} \boldsymbol{z}_1&\boldsymbol{z}_2\\\boldsymbol{z}_2^{\mathsf T}&z_3\end{pmatrix}\succeq 0,\ \boldsymbol{Y}\succeq 0$ with Constraints \eqref{eq:equal}--\eqref{eq:Sigma}, respectively. The Lagrangian function of \eqref{eq:inner_bound} has the following form:
\begin{align}
\label{Eq:lag}
&L(\boldsymbol{p},\boldsymbol{\tau},s, \boldsymbol{u}, \boldsymbol{Z},\boldsymbol{Y})=\sum_{k=1}^Kp_kQ_{t+1}(\boldsymbol{x}_t,\boldsymbol{\xi}_{t+1}^k)-s(\sum_{k=1}^Kp_k-1)+\boldsymbol{u}^{\mathsf T}(\boldsymbol{\tau}-\sum_{k=1}^K p_k\boldsymbol{\xi}_{t+1}^k)\nonumber\\
&+\begin{pmatrix}
    \boldsymbol{\Sigma}(\boldsymbol{x}_t) & \boldsymbol{\tau}-\boldsymbol{\mu}(\boldsymbol{x}_t)\\
    (\boldsymbol{\tau}-\boldsymbol{\mu}(\boldsymbol{x}_t))^{\mathsf T} & \gamma
    \end{pmatrix}\bullet \boldsymbol{Z} + \left(\eta \boldsymbol{\Sigma}(\boldsymbol{x}_t)-\sum_{k=1}^K p_k(\boldsymbol{\xi}_{t+1}^k-\boldsymbol{\mu}(\boldsymbol{x}_t))(\boldsymbol{\xi}_{t+1}^k-\boldsymbol{\mu}(\boldsymbol{x}_t))^{\mathsf T}\right)\bullet \boldsymbol{Y}\nonumber\\
& = \sum_{k=1}^K p_k\left(Q_{t+1}(\boldsymbol{x}_t,\boldsymbol{\xi}_{t+1}^k)-s-\boldsymbol{u}^{\mathsf T}\boldsymbol{\xi}_{t+1}^k-(\boldsymbol{\xi}_{t+1}^k-\boldsymbol{\mu}(\boldsymbol{x}_t))(\boldsymbol{\xi}_{t+1}^k-\boldsymbol{\mu}(\boldsymbol{x}_t))^{\mathsf T}\bullet \boldsymbol{Y}\right)\nonumber\\
&+\boldsymbol{\tau}^{\mathsf T}(\boldsymbol{u}+2\boldsymbol{z}_2)+s+\boldsymbol{\Sigma}(\boldsymbol{x}_t)\bullet \boldsymbol{z}_1 -\boldsymbol{\mu}(\boldsymbol{x}_t)^{\mathsf T}(2\boldsymbol{z}_2)+\gamma z_3 +\eta\boldsymbol{\Sigma}(\boldsymbol{x}_t)\bullet Y.
\end{align}
Because problem \eqref{eq:inner_bound} is convex and under the Slater's conditions, strong duality holds.
The maximization problem \eqref{eq:inner_bound} can be recast as
\begin{equation}
\label{eq:minmax}
    \min_{s,\boldsymbol{u},\boldsymbol{Z},\boldsymbol{Y}}\left\{\max_{\boldsymbol{p},\boldsymbol{\tau}}\left\{L(\boldsymbol{p},\boldsymbol{\tau},s, \boldsymbol{u}, \boldsymbol{Z},\boldsymbol{Y}): \boldsymbol{p\ge 0},\ \boldsymbol{\tau}\in\mathbb{R}^J\right\}\right\}.
\end{equation}
Following the Lagrangian function \eqref{Eq:lag}, after solving the inner maximization problem in \eqref{eq:minmax} over $\boldsymbol{p},\boldsymbol{\tau}$, we have
\begin{align*}
\min_{s,\boldsymbol{u},\boldsymbol{Z},\boldsymbol{Y}}\quad&s+\boldsymbol{\Sigma}(\boldsymbol{x}_t)\bullet \boldsymbol{z}_1 -\boldsymbol{\mu}(\boldsymbol{x}_t)^{\mathsf T}(2\boldsymbol{z}_2)+\gamma z_3 +\eta\boldsymbol{\Sigma}(\boldsymbol{x}_t)\bullet \boldsymbol{Y}\\
\text{s.t.}\quad& Q_{t+1}(\boldsymbol{x}_t,\boldsymbol{\xi}_{t+1}^k)-s-\boldsymbol{u}^{\mathsf T}\boldsymbol{\xi}_{t+1}^k-(\boldsymbol{\xi}_{t+1}^k-\boldsymbol{\mu}(\boldsymbol{x}_t))(\boldsymbol{\xi}_{t+1}^k-\boldsymbol{\mu}(\boldsymbol{x}_t))^{\mathsf T}\bullet \boldsymbol{Y}\le 0,\ \forall k\in[K],\\
& \boldsymbol{u}+2\boldsymbol{z}_2=0,\\
& \boldsymbol{Z}=\begin{pmatrix} \boldsymbol{z}_1&\boldsymbol{z}_2\\\boldsymbol{z}_2^{\mathsf T}&z_3\end{pmatrix}\succeq 0,\\
& \boldsymbol{Y}\succeq 0.
\end{align*}
Substituting $ \boldsymbol{u}=-2\boldsymbol{z}_2$ and combining with the outer minimization problem in \eqref{eq:bell_neutral}, we complete the proof. $\square$
\end{proof}

\begin{proof}[Theorem \ref{theorem:3-C}]
The proof follows the conic duality in functional spaces \citep{shapiro2001duality}. 
Using the ambiguity set defined in \eqref{eq:ambi_bound_C}, the inner maximization problem of \eqref{eq:bell_neutral} can be formulated as a conic linear program in a functional space as follows:
\begin{subequations}\label{eq:inner_3_C}
\begin{align}
\max_{P\in \mathcal{M}(\Xi_{t+1},\mathcal{F}_{t+1})}\quad&\mathbb{E}_P[Q_{t+1}(\boldsymbol{x}_t,\boldsymbol{\xi}_{t+1})]\\
\text{s.t.}\quad& \mathbb{E}_P[1]= 1,\label{eq:inner_3_C_1}\\
&(\mathbb{E}[\boldsymbol{\xi}_{t+1}]-\boldsymbol{\mu}(\boldsymbol{x}_t))^{\mathsf T}\boldsymbol{\Sigma}(\boldsymbol{x}_t)^{-1}(\mathbb{E}[\boldsymbol{\xi}_{t+1}]-\boldsymbol{\mu}(\boldsymbol{x}_t))\le \gamma,\\
& \mathbb{E}[(\boldsymbol{\xi}_{t+1}-\boldsymbol{\mu}(\boldsymbol{x}_t))(\boldsymbol{\xi}_{t+1}-\boldsymbol{\mu}(\boldsymbol{x}_t))^{\mathsf T}]\preceq \eta \boldsymbol{\Sigma}(\boldsymbol{x}_t)\label{eq:inner_3_C_3}
\end{align}
\end{subequations}
 We associate dual variables $s\in\mathbb{R}, \ \boldsymbol{Z}=\begin{pmatrix} \boldsymbol{z}_1&\boldsymbol{z}_2\\\boldsymbol{z}_2^{\mathsf T}&z_3\end{pmatrix}\succeq 0,\ \boldsymbol{Y}\succeq 0$ with Constraints \eqref{eq:inner_3_C_1}--\eqref{eq:inner_3_C_3}, respectively.
The Slater's constraint qualification conditions ensure that the primal problem has a non-empty relative interior. Therefore, strong duality holds and the dual problem can be formulated as \eqref{eq:bell_bound1_C}, which completes the proof. $\square$
\end{proof}

\end{document}